\documentclass[11pt,a4paper]{article}
\usepackage[T1]{fontenc}
\usepackage{mathptmx}
\usepackage{amsmath,amssymb,amsthm,mathtools}
\usepackage[margin=27mm,headheight=14pt,headsep=7mm]{geometry}
\usepackage{microtype}
\usepackage{fancyhdr}
\usepackage{bbm}
\usepackage{soul}
\usepackage{xcolor}   % »ò \usepackage{color}
\usepackage[colorlinks=true,linkcolor=blue,citecolor=blue,urlcolor=blue]{hyperref}
\makeatletter
\let\revision@texthl\hl
\newsavebox{\revision@mathbox}
\newcommand{\revision@mathhl}[2]{%
  \sbox{\revision@mathbox}{$\m@th#1#2$}%
  \begingroup\setlength{\fboxsep}{0pt}%
  \raisebox{0pt}[\ht\revision@mathbox][\dp\revision@mathbox]{%
    \colorbox{yellow}{\vphantom{\strut}\usebox{\revision@mathbox}}}%
  \endgroup}
\DeclareRobustCommand{\hl}[1]{%
  \ifmmode\mathpalette\revision@mathhl{#1}%
  \else\revision@texthl{#1}\fi}
\makeatother
\soulregister\cite7

\DeclareMathAlphabet{\mathcal}{OMS}{cmsy}{m}{n}
\hypersetup{pdftitle={Optimal orders for Hilbert-valued Littlewood--Paley functions of symmetric semigroups},pdfauthor={}}
\numberwithin{equation}{section}
\newtheorem{theorem}{Theorem}[section]
\newtheorem{corollary}[theorem]{Corollary}
\newtheorem{proposition}[theorem]{Proposition}
\newtheorem{lemma}[theorem]{Lemma}
\newtheorem{problem}[theorem]{Problem}
\theoremstyle{definition}

\theoremstyle{remark}
\newtheorem{remark}[theorem]{Remark}
\newcommand{\R}{\mathbb R}
\newcommand{\C}{\mathbb C}
\newcommand{\E}{\mathbb E}
\newcommand{\1}{\mathbf1}
\newcommand{\dd}{\,\mathrm d}
\newcommand{\norm}[2]{\left\|#1\right\|_{#2}}
\DeclareMathOperator*{\esssup}{ess\,sup}

\makeatletter
\renewcommand{\@seccntformat}[1]{\csname the#1\endcsname.\quad}
\makeatother
\allowdisplaybreaks

\begin{document}
\thispagestyle{empty}
\begin{flushleft}
\Large\bfseries Optimal orders for Hilbert-valued Littlewood--Paley functions of symmetric semigroups
\end{flushleft}

\noindent{Xudong Lai, Bo Li, Ji Li}
\smallskip

\noindent\textbf{Abstract.}
We study the Littlewood--Paley $g$-functions associated with symmetric
submarkovian semigroups and Hilbert-valued functions. For every
$1<p<\infty$ and $2\le q\le\infty$, we prove
\[
 \norm{G_q^T(f)}{L_p(\Omega)}
 \lesssim\max\{p',\sqrt p\}\norm{f-F(f)}{L_p(\Omega;\mathcal{H})},
\]
with the sharp order $\max\{p',\sqrt p\}$,
where the implicit constant is independent of $p,q$, the underlying space $\Omega$, the semigroup $(T_t)_{t>0}$,
and the Hilbert space $\mathcal{H}$.
Compactly supported radial functions show that, for each fixed $q>2$, the large-$p$ order cannot be replaced by $p^{1/q}$ in a dimension-free heat estimate.
In particular, for the case $q=2$, we answers affirmatively the question posed by Xu in his \emph{JFA} paper (2022) \cite[Problem~7]{XuJFA}.
For every finite $q>2$, our lower bound gives a negative answer to
the dimension-free question in Xu's \emph{JEMS} paper (2025)
\cite[Problem~1.8(i)]{XuJEMS}, already for real scalar functions.
Together, the upper and lower bounds establish the sharp
dimension-free order $\max\{p',\sqrt p\}$ for Hilbert-valued functions.

\par

\section{Introduction}

Let $(\Omega,\mathcal A,\mu)$ be a standard Borel $\sigma$-finite
measure space, and let $(T_t)_{t>0}$ be a symmetric submarkovian
semigroup. More precisely, $T_t=e^{-tA}$, where $A$ is a nonnegative
selfadjoint operator on $L_2(\Omega)$, $T_t$ is positive,
$T_t(\1)\le\1$, and $T_t$ has a measurable subprobability kernel for
every $t>0$. Its consistent extensions are contractions on
$L_p(\Omega)$ for $1\le p\le\infty$ and are strongly continuous for
$1\le p<\infty$. Throughout the paper, $\mathcal{H}$ is a real or complex
Hilbert space.
We use the same notation for the extension $T_t\otimes\mathrm{Id}_{\mathcal{H}}$ to $L_p(\Omega;\mathcal{H})$.

Let $F$ be the contractive projection onto the fixed point subspace
of $(T_t)_{t>0}$. For $2\le q\le \infty$, define the Littlewood--Paley $g$-function of $f$ by
\[
 G_q^T(f)(x)=
\begin{cases}
\displaystyle\left(\int_0^\infty
 \|t \partial_tT_t(f)(x)\|_{\mathcal{H}}^q
 \frac{\dd t}{t}\right)^{1/q},&2\le q<\infty, \\
\displaystyle\esssup_{t>0} \|t\partial_tT_t(f)(x)\|_{\mathcal{H}},&q=\infty.
\end{cases}
\]
The derivatives are the actual analytic semigroup derivatives on $L_p(\Omega;\mathcal{H})$.
As in \cite{XuJEMS}, let
$\mathrm L_{c,q,p}^T(\mathcal{H})$ be the least constant in
\begin{equation}\label{eq:best}
 \norm{G_q^T(f)}{L_p(\Omega)}
 \le\mathrm L_{c,q,p}^T(\mathcal{H})\norm{f-F(f)}{L_p(\Omega;\mathcal{H})},\quad 1<p<\infty.
\end{equation}

Stein's Littlewood--Paley--Stein theory \cite{Stein} established
$L_p$ bounds for $g$-functions associated with symmetric diffusion
semigroups. Subsequent developments connected square functions
with real-variable Hardy spaces \cite{FS}, probabilistic methods
\cite{Meyer}, and holomorphic functional calculus
\cite{Cowling,CDMY}. In the Banach-valued setting, one-sided
inequalities were characterized by martingale type and cotype,
first for Poisson semigroups and later for the original symmetric
diffusion semigroups \cite{MTX,Xu1998,XuII}. On complete Riemannian
manifolds, Coulhon, Duong, and Li \cite{CDL2003} obtained spatial
square-function estimates for $1<p\le2$ without geometric
assumptions. Further work addressed analytic regular contractions
\cite{LMX}, noncommutative $L_p$ spaces \cite{JLMX}, and optimal
dependence on $p$ \cite{XuJEMS,XuJFA}, with the dimension-free
Euclidean heat estimates raised as a problem in \cite{XuJEMS}.
Xu and Zhang \cite{XuZhangBG} used a comparison with martingale square functions to establish the optimal large-$p$ order for the reverse Littlewood--Paley--Stein inequality for symmetric diffusion semigroups. Hong, Xu, and Zhang \cite{HongXuZhang} obtained sharp estimates as $p\downarrow1$ for the reverse inequality for the classical Poisson semigroup, giving a partial answer to \cite[Problem~1.8(ii)]{XuJEMS}.
Recent endpoint advances include a geometry-independent weak
type $(1,1)$ bound for the full Poisson $g$-function
\cite{LLL2026} and related Riesz transform estimates in Euclidean
space \cite{OSS2026} and on complete noncompact Riemannian
manifolds \cite{CJLL2026}. For broader accounts and further
developments, including weighted inequalities and Banach-valued
extensions, see \cite{HNVW,SteinHA,Wilson} and the references therein.

In the scalar case $q=2$, the expected dimension-free order is
$\max\{p',\sqrt p\}$. This is also the question in
\cite[Problem~7]{XuJFA}. Our first result proves this estimate for
Hilbert-valued functions in the general symmetric submarkovian
setting and for every $q\ge2$, including the maximal derivative.

\begin{theorem}\label{thm:main}
Let $(T_t)_{t>0}$ be as above and let $\mathcal{H}$ be a Hilbert space. Then
\begin{equation}\label{eq:main}
 \mathrm L_{c,q,p}^T(\mathcal{H})\le
 \begin{cases}
 \displaystyle\frac12{(2^{11/2-6/q}p')}^{2-p},
       &1<p\le2,\\[4pt]
 \displaystyle\frac12{(2^{17/2-3/q}\sqrt p)}^{2-p'},
      &2\le p<\infty,
 \end{cases}
 \quad 2\le q\le\infty.
\end{equation}
Both bounds equal $1/2$ at $p=2$
and are independent of the underlying space $\Omega$, the semigroup $(T_t)_{t>0}$, and the Hilbert space $\mathcal{H}$.
In particular, the upper bounds are of order at most $p'$ as $p\downarrow1$ and $\sqrt p$ as $p\to\infty$,
uniformly in $q$.
\end{theorem}

For the heat semigroup, the large-$p$ order in Theorem~\ref{thm:main}
is necessary even when $q>2$. Write $H_t^{(d)}=e^{t\Delta_d}$ for
the Euclidean heat semigroup on $\R^d$.

\begin{theorem}\label{thm:sharp}
For every nonzero real or complex Hilbert space $\mathcal{H}$,
\begin{equation*}%
 \sup_{d\ge1}\mathrm L_{c,q,p}^{H^{(d)}}(\mathcal{H})
 \approx\max\{p',\sqrt p\},
 \quad 1<p<\infty,\quad 2\le q\le\infty.
\end{equation*}
The lower comparison constant may be taken to be $1/512$; the upper
bounds are those in \eqref{eq:main}. In the complex scalar case,
\begin{equation}\label{eq:explicitlower}
 \sup_{d\ge1}\mathrm L_{c,q,p}^{H^{(d)}}(\C)
 \ge\frac{\Gamma(q)^{1/q}}q\frac p{2\sqrt{p-1}},
 \quad 2\le q<\infty,
\end{equation}
while for $q=\infty$ the coefficient $\Gamma(q)^{1/q}/q$ is
replaced by $1/e$.
\end{theorem}

\begin{corollary}\label{cor:failure}
Fix $2<q\le\infty$. There is no constant
$C_q$ such that for every $1<p<\infty$, every $d\ge1$, and every $f\in L_p(\R^d)$,
\[
 \norm{G_q^{H^{(d)}}(f)}{L_p(\R^d)}
 \le C_q\max\{p',p^{1/q}\}\norm f{L_p(\R^d)}.
\]
This is true for real as well as complex scalar functions.
More precisely, for $p\ge2$,
\[
 \frac{\sup_{d\ge1}\mathrm L_{c,q,p}^{H^{(d)}}(\R)}
      {\max\{p',p^{1/q}\}}
 \ge\frac1{16e}p^{1/2-1/q}.
\]
\end{corollary}

For each finite $q>2$, this gives a negative answer to \cite[Problem~1.8(i)]{XuJEMS}
on the dimension-free forward estimate for the heat semigroup, already for $\mathcal{H}=\C$.
It does not exclude the estimate for a fixed Euclidean dimension.
The scalar and Hilbert-valued cases with $q=2$ follow from Theorem~\ref{thm:main}. The
subordinated Poisson estimate has a different large-$p$ behavior.

We briefly describe the proof. For small $1<p<2$, a change of jump rates realizes an exact time rescaling of a positive heat flow. Its first and second likelihood derivatives give two estimates in terms of entropy decreases. Integrating positive constant shifts converts these estimates into the required $L_p$ energy bounds. The second estimate controls the part of the maximal derivative exceeding the positive heat flow. This yields real-valued estimates for both $G_2^T$ and $G_\infty^T$ with order $p'$, and hence for all $G_q^T$ with $q\ge2$.

For large $p\ge 4$, we use the optimal scalar complex contraction angle. Independent signs determine rotations of successive time intervals and scalar phases formed from real-valued reverse martingale differences. Averaging the signs leaves the original dyadic derivatives and the conditional variances. One-dimensional Gaussian integration, together with conditional square-energy estimates for the martingale increments, yields an exponential estimate retaining the $L_2$ mass of the function. A quantitative interpolation argument yields the order $\sqrt p$. Positivity passes the dyadic bound to the maximal time derivative.

These estimates for small $1<p<2$ and large $p\ge 4$ are first proved for bounded conservative generators on probability spaces.
Complexification and interpolation in the space exponent give the real-valued estimates for the full range $1<p<\infty$. A single-atom conservative extension, finite-measure compression, and bounded Poisson approximation pass these estimates to the original semigroup. We then obtain the finite-dimensional real Hilbert-valued estimates by duality and Gaussian randomization, and approximation extends the estimates to arbitrary Hilbert spaces. Finally, we apply them to $f-F(f)$. Radial functions supported in thin Euclidean shells give the lower bound \eqref{eq:explicitlower}; a one-dimensional positive heat kernel estimate gives the order $p'$.

The remainder of this paper is organized as follows. Section~\ref{sec:prelim} contains the preliminaries. Sections~\ref{sec:low} and~\ref{sec:high} prove the real-valued estimates for smaller and larger $p$, respectively. Section~\ref{sec:all} obtains the real-valued estimates for all $1<p<\infty$. Section~\ref{sec:transfer} passes the estimates to general symmetric submarkovian semigroups and then to Hilbert-valued functions, proving Theorem~\ref{thm:main}. Section~\ref{sec:sharp} proves sharpness for the Euclidean heat semigroup, and Section~\ref{sec:applications} gives applications.

Throughout this paper, $L_p$ denotes a Bochner space when its value space is specified, and $F$ denotes the fixed point projection.
For simplicity, we write $p'=p/(p-1)$ and $u_t=T_t(f)$.
Unless otherwise specified, functions are understood to be scalar-valued;
for example, $L_p(\Omega)$ denotes the scalar-valued $L_p$ space,
whereas $L_p(\Omega;\mathcal{H})$ denotes the space of
$\mathcal{H}$-valued functions.
The relations $a\lesssim b$, $a\gtrsim b$,
and $a\approx b$ mean inequalities with absolute positive constants,
unless a dependence is explicitly indicated. Constants denoted by
$C$ may change from line to line. All time suprema are essential
suprema; on the bounded models below they agree with suprema of the
continuous representatives.

\section{Preliminaries}\label{sec:prelim}

\subsection{Analyticity, the fixed point projection, and the \texorpdfstring{$L_2$}{L2} identity}

For $1<p<\infty$, set $\theta_p=\arccos|1-2/p|$ and $\Sigma_{\theta_p}=\{z\in\C\setminus\{0\}:|\arg z|<\theta_p\}$.
The complex contraction theorem for symmetric contraction semigroups
\cite[Proposition~3]{CD} gives
\begin{equation}\label{eq:complex}
 \norm{T_z}{L_p(\Omega)\to L_p(\Omega)}\le1,
\quad z\in\Sigma_{\theta_p}.
\end{equation}
Only strictly interior angles are used below.
By the Marcinkiewicz--Zygmund extension theorem
\cite[Chapter~I, Section~8.24]{SteinHA},
\eqref{eq:complex} also holds on $L_p(\Omega;\mathcal{H})$
with the same constant; real Hilbert spaces are treated by
complexification.
The standard derivative estimates for bounded analytic semigroups
\cite[Chapter~2, Section~5]{Pazy} therefore give
\begin{equation}\label{eq:localderivatives}
 \norm{t^kA^kT_t}{L_p(\Omega;\mathcal{H})\to L_p(\Omega;\mathcal{H})}
 \le C_{p,k},\quad k\ge1.
\end{equation}
Here $C_{p,k}$ depends only on $p$ and $k$.
These bounds are used to identify limits and derivatives; their
dependence on $p$ is not used to prove \eqref{eq:main}.

Since $L_p(\Omega;\mathcal{H})$ is reflexive, the mean ergodic theorem gives
the contractive projection
\[
 F(f)=\lim_{s\to\infty}\frac1s\int_0^s T_t(f)\dd t
\]
onto the fixed point space. The projections are consistent on the
common domains of the $L_p$ realizations and agree with the spectral
projection of $A$ at $0$ on $L_2(\Omega;\mathcal{H})$. In particular,
$AT_t(F(f))=0$. Thus a bound with $\|f\|_{L_p(\Omega;\mathcal{H})}$ on the
right also holds with $\|f-F(f)\|_{L_p(\Omega;\mathcal{H})}$ there.
Hilbert-space spectral calculus and Tonelli's theorem give
\begin{equation}\label{eq:l2identity}
 \norm{G_2^T(f)}{L_2(\Omega)}^2
 =\frac14\norm{f-F(f)}{L_2(\Omega;\mathcal{H})}^2.
\end{equation}

We next collect several standard facts used below, starting with
the Gaussian moment identity \cite[Chapter~1]{Janson}.

\begin{lemma}\label{lem:gaussianmoments}
Let $\gamma_1,\ldots,\gamma_n$ be independent standard real Gaussian
random variables and $a_1,\ldots,a_n\in\R$. For every $p>0$,
\[
 \E\left|\sum_{k=1}^n a_k\gamma_k\right|^p
 = \left(\sum_{k=1}^n |a_k|^2\right)^{p/2}\E|\gamma_1|^p.
\]
The same identity holds for complex $a_k$ when the $\gamma_k$
are independent standard circular complex Gaussian random variables.
\end{lemma}

\begin{lemma}\label{lem:timeinterpolation}
For $2\le q\le\infty$ and almost every $x\in\Omega$,
\[
 G_q^T(f)(x)\le G_2^T(f)(x)^{2/q}
                    G_\infty^T(f)(x)^{1-2/q}.
\]
As a consequence, for $1<p<\infty$,
\[
\|G_q^T(f)\|_{L_p(\Omega)}
\le
\|G_2^T(f)\|_{L_p(\Omega)}^{2/q}
\|G_\infty^T(f)\|_{L_p(\Omega)}^{1-2/q}.
\]
\end{lemma}
\begin{proof}
This follows immediately from the definitions and H\"{o}lder's inequality.
\end{proof}

\begin{lemma}\label{lem:l2endpoint}
For every $2\le q\le\infty$,
\begin{equation}\label{eq:l2uniform}
 \|G_q^T(f)\|_{L_2(\Omega)}
 \le 2^{-5/4+1/(2q)}\|f-F(f)\|_{L_2(\Omega;\mathcal{H})}.
\end{equation}
For $q=2$, the inequality is an equality.
\end{lemma}
\begin{proof}
Spectral calculus, with $s=\log t$, gives
\[
 \int_{\R}\int_\Omega
 \|\partial_s (e^sAT_{e^s}(f))(x)\|_{\mathcal{H}}^2\dd\mu(x)\dd s
 =\frac18\|f-F(f)\|_{L_2(\Omega;\mathcal{H})}^2,
\]
because $\int_0^\infty (t\lambda)^2(1-t\lambda)^2e^{-2t\lambda}\dd t/t=1/8$ for $\lambda>0$, while the integral is zero at $\lambda=0$.
Together with \eqref{eq:l2identity}, this says that
$s\mapsto e^sAT_{e^s}(f)(x)$ belongs to $W^{1,2}(\R;\mathcal{H})$ for almost
every $x$, with limits zero at both ends. Integrating the derivative
of its squared norm on each side of $s$ yields
\[
 \sup_{s\in\R}\|e^sAT_{e^s}(f)(x)\|_{\mathcal{H}}^2
 \le\int_{\R}\|e^sAT_{e^s}(f)(x)\|_{\mathcal{H}}
       \|\partial_s  (e^sAT_{e^s}(f))(x)\|_{\mathcal{H}}\dd s,
\]
which, together with Cauchy--Schwarz's inequality, implies
$\|G_\infty^T(f)\|_{L_2(\Omega)}\le
2^{-5/4}\|f-F(f)\|_{L_2(\Omega;\mathcal{H})}$.
Lemma~\ref{lem:timeinterpolation} and
\eqref{eq:l2identity} now give \eqref{eq:l2uniform}.
The proof is complete.
\end{proof}

\subsection{Bounded conservative models}

Until Section~\ref{sec:transfer}, we work on a probability space
and assume
\begin{equation}\label{eq:model}
 A=\alpha(\mathrm{Id}_{\Omega}-P),\quad T_t=e^{-tA},\quad\alpha>0,
\end{equation}
where $P$ is a symmetric Markov kernel. This generator is bounded
on every $L_p(\Omega)$, including $L_\infty(\Omega)$. Its real-time kernels are
given by the Poisson series for $P$, and $T_z$ is entire as an
$L_\infty$ operator.
Since $A$ is bounded on $L_\infty(\Omega)$, the operator exponential
admits Taylor expansions in the $L_\infty$ operator norm. Thus, for each
fixed finite list of times and bounded functions, the Taylor remainders
are controlled in $L_\infty(\Omega)$.

For $0<t_1<\cdots<t_m$, consider the finite path law
\begin{equation*}%
 \dd\mu(x_0)T_{t_1}(x_0,\dd x_1)
 \prod_{j=2}^mT_{t_j-t_{j-1}}(x_{j-1},\dd x_j),
\end{equation*}
and denote its coordinates by $X_0,X_{t_1},\ldots,X_{t_m}$.
Reversibility and the Markov property imply
\begin{equation*}%
 \E[f(X_0)\mid X_{t_j},\ldots,X_{t_m}]
 =T_{t_j}(f)(X_{t_j}).
\end{equation*}
Read in reverse order, the right-hand side is a Hilbert-valued
martingale.

\begin{lemma}\label{lem:heatmax}
In the bounded conservative model, for $1<p<\infty$ and $g\ge0$,
\begin{equation*}%
 \norm{\sup_{t>0}T_t(g)}{L_p(\Omega)}
 \le p'\norm g{L_p(\Omega)}.
\end{equation*}
\end{lemma}
\begin{proof}
This is the standard maximal inequality obtained from Rota's
reverse-martingale representation \cite{Rota}; see also
\cite[Chapter~III]{Stein}.
\end{proof}

\section{Scalar estimates for smaller \texorpdfstring{$p$}{p}}\label{sec:low}

Throughout this section we work in the bounded conservative model
\eqref{eq:model} on a probability space.
The estimates below are independent of $\alpha$ and $P$.  Multiplying
the measure by a constant gives the same estimates on every finite
measure space.  The passage to general symmetric submarkovian
semigroups will be made in Section~\ref{sec:transfer}.

\subsection{Entropy dissipation and a change of jump rates}

For the moment, let $f$ be a bounded real-valued function with $f\geq\varepsilon>0$,
and put $ U_t^2=u_t^2+(tAu_t)^2$.
We use the entropy and the
temporal Fisher information
\[
 E(t)=\int_\Omega u_t(x)\log u_t(x)\dd \mu(x),
 \quad I(t)=\int_\Omega\frac{(Au_t(x))^2}{u_t(x)}\dd \mu(x).
\]
The Poisson expansion of $T_t$ gives representatives of $u$, $Au$, and
$A^2u$ that are smooth in $t>0$ outside one null set.
Define a measure \(m\) on ordered pairs by
\[
\mathrm dm(x,y)=\alpha\,\mathrm d\mu(x)P(x,\mathrm dy).
\]
The symmetry of \(P\) with respect to \(\mu\) implies that \(m\) is invariant under exchanging \(x\) and \(y\).
With $w_t=-Au_t/u_t$, the Markov property of $P$ and symmetry yield
\begin{align*}
 -E'(t)
 &=\frac12\int_{\Omega\times\Omega} (u_t(x)-u_t(y))(\log u_t(x)-\log u_t(y)) \dd m(x,y)\ge0,\\
 -I'(t)
 &=\frac12\int_{\Omega\times\Omega} (u_t(x)+u_t(y))(w_t(x)-w_t(y))^2 \dd m(x,y)\ge0.
\end{align*}

Fix $a>0$ and $\tau>0$.
For $|\theta|<1/2$ and $0\leq s\leq\tau$, define
\begin{align*}
 v_\theta(s,x)&=u_{a+(1+\theta)s}(x),\\
 h_\theta(s,x,y)&=\frac{v_\theta(s,x)-v_\theta(s,y)}
 {v_\theta(s,x)+v_\theta(s,y)},\\
 q_\theta(s,x,\dd y)&=\alpha  (1+\theta h_\theta(s,x,y))P(x,\dd y).
\end{align*}
These are genuine nonnegative jump rates, with total rate at most
$3\alpha/2$.  Start the process with density $u_a/\|f\|_{L_1(\Omega)}$ relative to $\mu$.
At density $v_\theta/\|f\|_{L_1(\Omega)}$, the net probability
flux (incoming minus outgoing) is a signed measure in $x$.
Using the symmetry of $m$ and the antisymmetry of $h_\theta$,
we obtain
\begin{align*}
&\int_\Omega
 \left[ \frac{v_\theta(s,y)}{\|f\|_{L_1(\Omega)}}\dd \mu(y) \right]q_\theta(s,y,\dd x)
  - \int_\Omega \left[ \frac{v_\theta(s,x)}{\|f\|_{L_1(\Omega)}}\dd \mu(x) \right]q_\theta(s,x,\dd y)\\
&\quad= \frac{\alpha\dd \mu(x)}{\|f\|_{L_1(\Omega)}} \int_\Omega [ v_\theta(s,y)(1+\theta h_\theta(s,y,x))
 - v_\theta(s,x)(1+\theta h_\theta(s,x,y)) ]P(x,\dd y)\\
&\quad= \frac{\alpha(1+\theta)\dd \mu(x)}  {\|f\|_{L_1(\Omega)}} \int_\Omega (v_\theta(s,y)-v_\theta(s,x))P(x,\dd y)\\
&\quad= \frac{\alpha(1+\theta)}{\|f\|_{L_1(\Omega)}}
 (P-\mathrm{Id}_{\Omega})v_\theta(s,x)\dd \mu(x).
\end{align*}
Thus the net flux density with respect to $\mu$ is
$-(1+\theta)Av_\theta/
\|f\|_{L_1(\Omega)}$.
Uniqueness for the bounded nonautonomous forward equation therefore
gives the marginal density $v_\theta(s)/\|f\|_{L_1(\Omega)}$ at every time $s$.
In the standard Borel setting, the process can be constructed using a
Poisson clock of rate $3\alpha/2$ and rejection.  We retain self-jump
marks; they have $h_\theta=0$ and make no contribution to the expressions
below.

Let $(\xi_s)_{0\leq s\leq\tau}$ denote the process for $\theta=0$,
and let $N(\dd s,\dd y)$ be its marked jump measure.  Its compensator is
$q_0(s,\xi_{s-},\dd y)\dd s=\alpha P(\xi_{s-},\dd y)\dd s$.
The likelihood of the perturbed marked path relative to this path law is
\[
 R_\theta=
 \exp\left(-\int_0^\tau\!\int_\Omega
 \theta h_\theta(s,\xi_{s-},y)q_0(s,\xi_{s-},\dd y)\dd s\right)
 \prod_{(s,y)\in N}(1+\theta h_\theta(s,\xi_{s-},y)).
\]
Differentiating $\log R_\theta$ with respect to $\theta$ at $\theta=0$ gives
\[
R'_0=V, \quad R''_0=V^2-[V]+2W.
\]
Here $V$ and $W$ are compensated stochastic integrals,
and $[V]$ and $[W]$ denote the optional quadratic variations at time $\tau$ of the martingales defining $V$ and $W$, respectively,
given explicitly by
\[
\begin{aligned}
V&=\int_0^\tau\!\int_\Omega h(s,\xi_{s-},y)[N(\dd s,\dd y)-q_0(s,\xi_{s-},\dd y)\dd s],\quad
 [V]= \int_0^\tau\!\int_\Omega h(s,\xi_{s-},y)^2N(\dd s,\dd y),\\
W&=\int_0^\tau\!\int_\Omega\dot h(s,\xi_{s-},y)[N(\dd s,\dd y)-q_0(s,\xi_{s-},\dd y)\dd s],\quad
 [W]=\int_0^\tau\!\int_\Omega \dot h(s,\xi_{s-},y)^2\,N(\mathrm ds,\mathrm dy),
\end{aligned}
\]
where $h=h_0$ and $\dot h=\partial_\theta h_\theta|_{\theta=0}$.
All integrands are predictable, with the state variable evaluated
at the pre-jump position $\xi_{s-}$.
The use of the optional quadratic variation is
essential in this pathwise identity.
On this fixed interval, the positive lower bound for $u$, the boundedness
of $A$, and the exponential moments of the number of marks provide an
integrable bound for the first two likelihood derivatives, locally
uniformly in $\theta$.  Hence these derivatives may also be taken under
expectations and in the endpoint marginal.

Martingale isometry and symmetry give
\begin{align}
 \mathbb E V^2
 &=\mathbb E([V])
 =\frac{1}{2 \|f\|_{L_1(\Omega)}}\int_0^\tau\!\int_{\Omega\times\Omega}
 \frac{(u_{a+s}(x)-u_{a+s}(y))^2}
 {u_{a+s}(x)+u_{a+s}(y)}\dd m(x,y) \dd s
 \leq\frac{E(a)-E(a+\tau)}{2 \|f\|_{L_1(\Omega)}},\label{low:scorefirst}\\
 \mathbb E W^2
 &=\mathbb E([W])=\frac{2}{\|f\|_{L_1(\Omega)}}\int_0^\tau s^2\!\int_{\Omega\times\Omega} \frac{u_{a+s}(x)^2u_{a+s}(y)^2} {(u_{a+s}(x)+u_{a+s}(y))^3}
 (w_{a+s}(x)-w_{a+s}(y))^2\dd m(x,y)\dd s\notag\\
 &\leq \frac{1}{4\|f\|_{L_1(\Omega)}}\int_0^\tau s^2[-I'(a+s)] \dd s
 \leq \frac{\tau^2 I(a)}{4\|f\|_{L_1(\Omega)}}.\notag
\end{align}
The first comparison follows from $(b-d)(\log b-\log d)\geq2(b-d)^2/(b+d)$ for $b,d>0$,
while the second uses $16b^2d^2\leq(b+d)^4$.
All factors of $\|f\|_{L_1(\Omega)}$ arise from the probability
normalization of the initial density.

\subsection{Two adjacent time intervals}

In this subsection we combine the likelihood identities with entropy dissipation on the adjacent intervals $[t/4,t/2]$ and $[t/2,t]$.
This gives a joint estimate for the first and second time derivatives and a temporal Fisher information estimate.

\begin{lemma}\label{lem:low:twointervals}
Let $f$ be a bounded real-valued function with $f\geq\varepsilon>0$.
For every $t>0$ we have
\[\int_\Omega \frac{|tAu_t(x)|^2}{u_t(x)}\dd \mu(x)\leq 2[E(t/2)-E(t)],\quad
\int_\Omega \frac{|tAu_t(x)\,t^2A^2u_t(x)|}{U_t(x)}\dd \mu(x)
\leq 8[E(t/4)-E(t)].
\]
\end{lemma}

\begin{proof}
Fix $t>0$ and invoke the preceding construction on $[t/2,t]$,
with $a=\tau=t/2$.
Let $Y=\xi_{t/2}$.
Differentiation of the endpoint marginal, whose density is  $u_{t+\theta t/2}/\|f\|_{L_1(\Omega)}$, now yields
\begin{equation}\label{low:endpointfirst}
\mathbb E(V\mid Y)
=-\frac{tAu_t(Y)}{2u_t(Y)},\quad
\mathbb E(V^2-[V]+2W\mid Y)
=\frac{t^2A^2u_t(Y)}{4u_t(Y)}.
\end{equation}
These identities are understood as identities of Radon--Nikodym derivatives,
so no choice of conditional probability kernels is needed.

Since $Y$ has density $u_t/\|f\|_{L_1(\Omega)}$, the conditional Cauchy--Schwarz inequality and \eqref{low:scorefirst} give
\[
\frac{t^2}{4} I(t)
=\int_\Omega \frac{|tAu_t(x)|^2}{4u_t(x)}\dd \mu(x)
=\|f\|_{L_1(\Omega)}\mathbb E\!\left[\left|\frac{tAu_t(Y)}{2u_t(Y)}\right|^2\right]
\le \|f\|_{L_1(\Omega)}\mathbb E V^2
\leq \frac{E(t/2)-E(t)}{2},
\]
which  proves the first desired conclusion.
The same estimate with $t$ replaced by $t/2$ gives
\[
\left(\frac{t}{2}\right)^2 I(t/2)
\leq 2[E(t/4)-E(t/2)].
\]
Consequently,
the moment estimates from the preceding subsection imply
\[
\mathbb E V^2 \leq \frac{E(t/2)-E(t)}{2\|f\|_{L_1(\Omega)}},\quad
\mathbb E W^2 \leq \frac{t^2 I(t/2)}{16\|f\|_{L_1(\Omega)}} \leq \frac{E(t/4)-E(t/2)}{2\|f\|_{L_1(\Omega)}}.
\]

To prove the second estimate, we deduce from \eqref{low:endpointfirst}, the conditional Jensen inequality and the Cauchy--Schwarz inequality that
\begin{align*}
\int_\Omega \frac{|tAu_t(x)\,t^2A^2u_t(x)|}{U_t(x)}\,\mathrm d\mu(x)
&= \|f\|_{L_1(\Omega)}\mathbb E\!\left[\frac{\left|\dfrac{tAu_t(Y)}{u_t(Y)}\dfrac{t^2A^2u_t(Y)}{u_t(Y)}\right|}{\sqrt{1+\left(\frac{tAu_t(Y)}{u_t(Y)}\right)^2}}\right] \notag\\
&= 4\|f\|_{L_1(\Omega)}\mathbb E\!\left[\sqrt{\frac{4(\mathbb E(V\mid Y))^2}{1+4(\mathbb E(V\mid Y))^2}}|\mathbb E(V^2-[V]+2W\mid Y)|\right] \notag\\
&\leq 4\|f\|_{L_1(\Omega)}\mathbb E\!\left[\sqrt{\frac{4(\mathbb E(V\mid Y))^2}{1+4(\mathbb E(V\mid Y))^2}}|\mathbb E(V^2-[V]\mid Y)|\right]\notag \\
&\quad +8\|f\|_{L_1(\Omega)}\mathbb E\!\left[\sqrt{\frac{4(\mathbb E(V\mid Y))^2}{1+4(\mathbb E(V\mid Y))^2}}|\mathbb E(W\mid Y)|\right]\notag \\
&\leq 4\|f\|_{L_1(\Omega)}\mathbb E\![\mathbb E(V^2+[V]\mid Y)]+16\|f\|_{L_1(\Omega)}\mathbb E\![|\mathbb E(V\mid Y)|\,|\mathbb E(W\mid Y)|]\notag \\
&\leq 8\|f\|_{L_1(\Omega)}\mathbb E V^2+16\|f\|_{L_1(\Omega)}(\mathbb E V^2)^{1/2}(\mathbb E W^2)^{1/2}\notag \\
&\le 8[E(t/4)-E(t)].
\end{align*}
The last two inequalities use the arithmetic--geometric mean inequality and the monotonicity of $E$.
No independence or orthogonality between $V$ and $W$ is used.
The proof is complete.
\end{proof}

\subsection{Positive shifts and integral estimates}

We apply the preceding entropy estimates to positive shifts of the initial datum and integrate over the shift parameter.
This converts entropy differences into differences of $L_p$ norms and gives the following space--time energy estimates.

\begin{proposition}\label{prop:low:energies}
For every $1<p<2$, if real-valued $f\geq0$ is bounded, then
\begin{align*}
M_p&=\int_0^\infty\!\int_\Omega u_t(x)^{p-2}|tAu_t(x)|^2\dd \mu(x)\frac{\dd t}{t}
\leq\frac{\log4}{p(p-1)}\|f\|_{L_p(\Omega)}^p, \\
N_p&=\int_0^\infty\!\int_\Omega U_t(x)^{p-2}|tAu_t(x)\,t^2A^2u_t(x)|\dd \mu(x)\frac{\dd t}{t}
\leq\frac{8\log4}{p(p-1)}\|f\|_{L_p(\Omega)}^p.
\end{align*}
\end{proposition}

\begin{proof}
First suppose $f\geq\varepsilon>0$.
For $s>0$,
conservativity gives $T_t(f+s)=u_t+s$,
whereas $tAu_t$ and $t^2A^2u_t$ are unchanged.
Define the nonnegative remainder
\[
\psi_s(v)=\int_0^v\frac{v-z}{z+s}\dd z=(v+s)\log(1+v/s)-v, \quad v\geq0.
\]
For $v>0$, Tonelli's theorem and Euler's beta integral yield
\[
\int_0^\infty s^{p-2}\psi_s(v) \dd s
=\int_0^v(v-z)\left(\int_0^\infty\frac{s^{p-2}}{z+s} \dd s\right)\dd z
=\frac{-\pi}{p(p-1)\sin(\pi p)}v^p.
\]
The same identity holds at $v=0$.
Moreover, set
\[
E_s(r)=\int_\Omega(u_r(x)+s)\log(u_r(x)+s)\dd \mu(x).
\]
Since $(u_r+s)\log(u_r+s)=\psi_s(u_r)+(u_r+s)\log s+u_r$, conservativity of $T_t$ gives
\[
E_s(a)-E_s(b)
=\int_\Omega\psi_s(u_a(x))\dd \mu(x)-\int_\Omega\psi_s(u_b(x))\dd \mu(x).
\]
Integrating the two nonnegative remainders separately by Tonelli's theorem and then subtracting,
we obtain
\begin{equation}
\int_0^\infty s^{p-2}[E_s(a)-E_s(b)] \dd s
=\frac{-\pi}{p(p-1)\sin(\pi p)}
\left[\|u_a\|_{L_p(\Omega)}^p-\|u_b\|_{L_p(\Omega)}^p\right],
\quad 0<a<b.
\label{low:entropybeta}
\end{equation}

Apply Lemma \ref{lem:low:twointervals} to $f+s$.
Using $(u_t(x)+s)^2+(tAu_t(x))^2\leq (U_t(x)+s)^2$,
Tonelli's theorem and \eqref{low:entropybeta} give
\begin{align*}
\int_\Omega U_t(x)^{p-2}|tAu_t(x)\,t^2A^2u_t(x)|\dd \mu(x)
&=-\frac{\sin(\pi p)}{\pi}\int_0^\infty s^{p-2}\int_\Omega\frac{|tAu_t(x)\,t^2A^2u_t(x)|}{U_t(x)+s}\dd \mu(x) \dd s\\
&\leq -\frac{\sin(\pi p)}{\pi}\int_0^\infty s^{p-2}\int_\Omega\frac{|tAu_t(x)\,t^2A^2u_t(x)|}{\sqrt{(u_t(x)+s)^2+(tAu_t(x))^2}}\dd \mu(x) \dd s\\
&\leq-\frac{8\sin(\pi p)}{\pi}\int_0^\infty s^{p-2}[E_s(t/4)-E_s(t)] \dd s\\
&=\frac{8}{p(p-1)}\left[\|u_{t/4}\|_{L_p(\Omega)}^p-\|u_t\|_{L_p(\Omega)}^p\right].
\end{align*}
Likewise, we also obtain
\[
\int_\Omega u_t(x)^{p-2}|tAu_t(x)|^2\dd \mu(x)\le \frac{2}{p(p-1)}\left[\|u_{t/2}\|_{L_p(\Omega)}^p-\|u_t\|_{L_p(\Omega)}^p\right].
\]

The semigroup property and $L_p$ contractivity imply that $t\mapsto\|u_t\|_{L_p(\Omega)}^p$ is decreasing,
and strong continuity gives $\lim_{t\downarrow0}\|u_t\|_{L_p(\Omega)}^p=\|f\|_{L_p(\Omega)}^p$.
For any $K>1$, Tonelli's theorem applied to the nonnegative Stieltjes measure $-\dd\|u_s\|_{L_p(\Omega)}^p$ gives
\[
\int_0^\infty\left[\|u_{t/K}\|_{L_p(\Omega)}^p-\|u_t\|_{L_p(\Omega)}^p\right]\frac{\dd t}{t}
=\int_{(0,\infty)}\left(\int_s^{Ks}\frac{\dd t}{t}\right)\left(-\dd\|u_s\|_{L_p(\Omega)}^p\right)
\leq\log K\,\|f\|_{L_p(\Omega)}^p.
\]
Integrating the two preceding estimates against $\dd t/t$ and taking $K=4$ and $K=2$,
respectively, proves this proposition provided $f\geq\varepsilon>0$.

Finally,
let $f\geq0$ be bounded.
Apply the estimates to $f+\varepsilon$ and let $\varepsilon\downarrow0$.
At $U_t(x)=0$,
define the integrand in $N_p$ to be zero.
This is a continuous extension because
$|tAu_t(x)|U_t(x)^{p-2}|t^2A^2u_t(x)|
\leq U_t(x)^{p-1}|t^2A^2u_t(x)|$.
For the other  integrand in  $M_p$, at points where $u_t(x)=0$, use the convention
\[
u_t(x)^{p-2}|tAu_t(x)|^2=
\begin{cases}
0,&tAu_t(x)=0,\\
+\infty,& tAu_t(x)\ne0.
\end{cases}
\]
This extension is lower semicontinuous.
Since $\|f+\varepsilon\|_{L_p(\Omega)}^p\to\|f\|_{L_p(\Omega)}^p$,
Fatou's lemma gives both estimates for bounded $f\geq0$.
The proof is complete.
\end{proof}

\subsection{Estimates for real-valued data}

Combining the heat maximal estimate with Proposition~\ref{prop:low:energies},
we first estimate the square function and the derivative maximum for nonnegative data.
Interpolation and decomposition into positive and negative parts then give the following estimate.

\begin{proposition}\label{prop:low:scalar}
For every $1<p<2$, if $f\in L_p(\Omega)$ is real-valued, then
\[
\|G_q^T(f)\|_{L_p(\Omega)}\leq2^{7/2-6/q}p'\|f\|_{L_p(\Omega)},
\quad 2\leq q\leq\infty.
\]
\end{proposition}

\begin{proof}
\noindent\textbf{Step 1: Estimates for nonnegative bounded data.} We begin with bounded $f\geq0$.

\noindent\textbf{Substep 1.1: The square function bound.} The heat maximal estimate and the energy estimate in Proposition~\ref{prop:low:energies} control the square function,
namely,
\[
G_2^T(f)(x)^2
\leq\left(\sup_{t>0}u_t(x)\right)^{2-p}\int_0^\infty u_t(x)^{p-2}|tAu_t(x)|^2\frac{\dd t}{t}.
\]
H\"older's inequality in $\Omega$ with exponents $2/(2-p)$ and $2/p$,
followed by Proposition~\ref{prop:low:energies},
yields
\begin{equation}
\|G_2^T(f)\|_{L_p(\Omega)}
\leq(p')^{1-p/2}\left(\frac{\log4}{p(p-1)}\right)^{1/2}\|f\|_{L_p(\Omega)}.
\label{low:square}
\end{equation}

\noindent\textbf{Substep 1.2: The derivative maximal bound.}
In the logarithmic variable $r=\log t$,
set
\[
v(r,x)=(|e^rAu_{e^r}(x)|-u_{e^r}(x))_+.
\]
Since $\partial_r (u_{e^r})=-e^rAu_{e^r}$ and $\partial_r(e^rAu_{e^r})=e^rAu_{e^r}-e^{2r}A^2u_{e^r}$,
the chain rule for locally absolutely continuous scalar functions gives,
almost everywhere,
\begin{align*}
|\partial_r(v(r,x)^p)|
&\leq p|e^rAu_{e^r}(x)|^{p-1}(2|e^rAu_{e^r}(x)|+|e^{2r}A^2u_{e^r}(x)|)\mathbbm{1}_{\{|e^rAu_{e^r}|>u_{e^r}\}}\\
&\leq p\left(2u_{e^r}(x)^{p-2}|e^rAu_{e^r}(x)|^2+2^{(2-p)/2}U_{e^r}(x)^{p-2}|e^rAu_{e^r}(x)e^{2r}A^2u_{e^r}(x)|\right)\mathbbm{1}_{\{|e^rAu_{e^r}|>u_{e^r}\}},
\end{align*}
where the zero-value conventions above apply to the first term.

We justify the endpoints in this variation estimate.
Since $f$ is bounded and $\mu$ is finite,
the spectral theorem gives
\begin{align*}
\int_0^\infty\|tAu_t\|_{L_2(\Omega)}^2\frac{\dd t}{t}
&=\frac14\|f-F(f)\|_{L_2(\Omega)}^2,\\
\int_0^\infty\|t^2A^2u_t\|_{L_2(\Omega)}^2\frac{\dd t}{t}
&=\frac38\|f-F(f)\|_{L_2(\Omega)}^2,\\
\int_{\mathbb R}\|\partial_r(e^rAu_{e^r})\|_{L_2(\Omega)}^2 \dd r
&=\frac18\|f-F(f)\|_{L_2(\Omega)}^2.
\end{align*}
Thus,
for almost every $x$,
$e^rAu_{e^r}(x)$ belongs to $W^{1,2}(\mathbb R)$ and its continuous representative tends to zero at both ends.
The function $v(\cdot,x)$ is locally absolutely continuous,
satisfies $0\leq v(r,x)\leq|e^rAu_{e^r}(x)|$,
and has derivative bounded by $2|e^rAu_{e^r}(x)|+|e^{2r}A^2u_{e^r}(x)|$.
Hence it too belongs to $W^{1,2}(\mathbb R)$ and has zero limits at both ends.
The total variation of the nonnegative function $v(\cdot,x)^p$ satisfies
\begin{align*}
\int_{\mathbb R}|\partial_r(v(r,x)^p)| \dd r
&=\int_{-\infty}^{\ell}|\partial_r(v(r,x)^p)| \dd r+\int_\ell^\infty|\partial_r(v(r,x)^p)| \dd r\\
&\geq\left|\int_{-\infty}^{\ell}\partial_r(v(r,x)^p) \dd r\right|+\left|\int_\ell^\infty\partial_r(v(r,x)^p) \dd r\right|=2v(\ell,x)^p,
\end{align*}
for all $\ell\in\mathbb R$.
This total variation estimate,
combined with the bound on $|\partial_r(v(r,x)^p)|$ above and Proposition~\ref{prop:low:energies},
implies
\begin{align}\label{low:excess}
\left\|\sup_{r\in\mathbb R}v(r,\cdot)\right\|_{L_p(\Omega)}^p\leq pM_p+\frac p2\,2^{(2-p)/2}N_p\leq\frac{2\log2(1+4\sqrt2)}{p-1}\|f\|_{L_p(\Omega)}^p.
\end{align}
Note that
\[
G_\infty^T(f)(x)\leq\sup_{t>0}u_t(x)+\sup_{r\in\mathbb R}v(r,x).
\]
The heat maximal estimate in Lemma~\ref{lem:heatmax},
together with \eqref{low:excess},
consequently gives
\begin{equation}
\|G_\infty^T(f)\|_{L_p(\Omega)}
\leq\left[p'+\left(\frac{2\log2(1+4\sqrt2)}{p-1}\right)^{1/p}\right]\|f\|_{L_p(\Omega)}.
\label{low:maximum}
\end{equation}

\noindent\textbf{Substep 1.3: Simplification of the endpoint constants.}
We  simplify the constants in \eqref{low:square} and \eqref{low:maximum}.

For the coefficient in \eqref{low:square},
the function $\log p-(1-1/p)\log2$ is increasing and is zero at $p=1$,
so $p^{-1}\leq2^{1/p-1}$.
Since $0<(p-1)/p<1$ and $\log4<2$, one has
\[
(p')^{1-p/2}\left(\frac{\log4}{p(p-1)}\right)^{1/2}
=\frac{\sqrt{\log4}}p\left(\frac{p-1}p\right)^{(p-1)/2}p'
\leq\frac{\sqrt2}{p}p'\leq2^{1/p-1/2}p'.
\]
Thus,
for bounded nonnegative $f$, we have
\[
\|G_2^T(f)\|_{L_p(\Omega)}
\leq2^{1/p-1/2}p'\|f\|_{L_p(\Omega)}.
\]

For the coefficient in \eqref{low:maximum},
note that $(p-1)^{-1/p}\leq(p-1)^{-1}\leq p'$ and $1+2\log2(1+4\sqrt2)<8\sqrt2$.
H\"older's inequality for a two-term sum therefore gives
\begin{align*}
p'+\left(\frac{2\log2(1+4\sqrt2)}{p-1}\right)^{1/p}
&\leq\left(1+[2\log2(1+4\sqrt2)]^{1/p}\right)p'\\
&\leq2^{1-1/p}\left(1+2\log2(1+4\sqrt2)\right)^{1/p}p'\leq2^{1/p+5/2}p'.
\end{align*}
Consequently, one has
\[
\|G_\infty^T(f)\|_{L_p(\Omega)}
\leq2^{1/p+5/2}p'\|f\|_{L_p(\Omega)}.
\]

\noindent\textbf{Substep 1.4: Interpolation between the endpoint bounds.}
For $2<q<\infty$,
Lemma~\ref{lem:timeinterpolation},
followed by the simplified endpoint estimates in Substep 1.3,
yields
\[
\|G_q^T(f)\|_{L_p(\Omega)}
\leq\|G_2^T(f)\|_{L_p(\Omega)}^{2/q}\|G_\infty^T(f)\|_{L_p(\Omega)}^{1-2/q}\leq2^{1/p+5/2-6/q}p'\|f\|_{L_p(\Omega)}.
\]
The same bound holds at $q=2$ and $q=\infty$ by the endpoint estimates above.

\noindent\textbf{Step 2: Extension to real $L_p$ data.}
For bounded real $f$,
write $f=f_+-f_-$.
The triangle inequality and H\"older's inequality give
\[
\|G_q^T(f)\|_{L_p(\Omega)}
\leq\|G_q^T(f_+)\|_{L_p(\Omega)}+\|G_q^T(f_-)\|_{L_p(\Omega)},\quad
\|f_+\|_{L_p(\Omega)}+\|f_-\|_{L_p(\Omega)}
\leq2^{1-1/p}\|f\|_{L_p(\Omega)}.
\]
Applying the nonnegative estimate from Substep 1.4 to $f_+$ and $f_-$ gives
\begin{align*}
\|G_q^T(f)\|_{L_p(\Omega)}
&\leq2^{1/p+5/2-6/q}p'\left(\|f_+\|_{L_p(\Omega)}+\|f_-\|_{L_p(\Omega)}\right)\\
&\leq2^{1/p+5/2-6/q}2^{1-1/p}p'\|f\|_{L_p(\Omega)}=2^{7/2-6/q}p'\|f\|_{L_p(\Omega)}.
\end{align*}
This proves the estimate for bounded real-valued $f$.

To complete the proof,
we extend the estimates from bounded real functions to arbitrary real $f\in L_p(\Omega)$ by approximation.
Indeed, $A$ is bounded on $L_p$,
so the actual derivatives converge on every compact positive time interval,
and Fatou's lemma applies after passage to a subsequence.
\end{proof}

\section{Scalar estimates for larger \texorpdfstring{\(p\)}{p}}
\label{sec:high}

We continue with the bounded conservative model \eqref{eq:model} and consider real-valued functions throughout this section. All constants obtained below are independent of $\alpha$ and $P$, and of the number and location of the time scales. For $s>0$ and a positive integer $N$, put
\[
B_N(f)=\left(\sum_{j=1}^N|2^jsAu_{2^{j+1}s}|^2\right)^{1/2}.
\]
We will first prove an exponential estimate for $B_N(f)$ which retains the $L_2$ mass of $f$. This is the point needed for the subsequent interpolation.

\subsection{A finite-path phase estimate}

Use the times $2^js$, $1\le j\le N+1$, and consider the probability law on $\Omega^{N+2}$ given by
\begin{equation}\label{high:path-law}
\dd\mu(\omega_0)T_{2s}(\omega_0,\dd\omega_1)\prod_{j=2}^{N+1}T_{2^{j-1}s}(\omega_{j-1},\dd\omega_j).
\end{equation}
Each marginal is $\mu$. Detailed balance and the Markov property imply
\[
\E[f(\omega_0)\mid\omega_j,\ldots,\omega_{N+1}]=u_{2^js}(\omega_j).
\]
Consequently $u_{2^{N+1}s}(\omega_{N+1}),\ldots,u_{2s}(\omega_1)$, in that order, is a finite real-valued martingale. Write
\[
D_j=u_{2^js}(\omega_j)-u_{2^{j+1}s}(\omega_{j+1}),\quad 1\le j\le N.
\]
We denote by $\E_\omega$ expectation in the same finite path with initial point $\omega_0=\omega$, and define the nonnegative functions $C_j$ by
\[
C_j(\omega)=\E_\omega[(D_j-\E_\omega D_j)^2],\quad 1\le j\le N.
\]

\begin{lemma}\label{high:phase}
For every bounded real-valued $f$ and every $\gamma_1,\ldots,\gamma_N\in\R$, one has
\begin{equation}\label{high:phase-bound}
\int_\Omega\exp\left[\sum_{j=1}^N\left(\frac{\gamma_j}{3}2^{j+1}sAu_{2^{j+1}s}(\omega)-\frac12C_j(\omega)\gamma_j^2\right)\right]\dd\mu(\omega)\le1.
\end{equation}
\end{lemma}

\begin{proof}
\noindent\textbf{Step 1: Construction and contractivity of $U_R$.}
Fix $0<\delta<1/3$ and $N+1$ independent Rademacher signs $\varepsilon_1,\ldots,\varepsilon_{N+1}$. For the moment their values are fixed. Set
\[
\alpha_j=
\begin{cases}\delta\varepsilon_1,&j=1,\\
\delta(2\varepsilon_j-\varepsilon_{j-1}),&2\le j\le N+1.
\end{cases}
\]
In particular $|\alpha_j|\le3\delta<1$ for every $j$. For $R>2$ and $j\ge2$ put
\[
\theta_R=\arccos|1-2/R|=2\arcsin R^{-1/2},\quad z_1=2s e^{i\alpha_1\theta_R},\quad z_j=2^{j-1}s e^{i\alpha_j\theta_R}.
\]
The scalar complex-time contraction theorem \cite[Proposition~3]{CD} gives $\|T_{z_j}\|_{L_R(\Omega)\to L_R(\Omega)}\le1$: all the angles used here are strictly inside the contraction sector. With the conventions $\varepsilon_0=1$ and $\gamma_0=\gamma_{N+1}=0$, define
\[
\kappa_j=\varepsilon_j\gamma_j-\varepsilon_{j-1}\gamma_{j-1},\quad 1\le j\le N+1.
\]
The functions $\beta_j=\kappa_j u_{2^js}$ are real and bounded. Let
\begin{equation}\label{high:nested-phase}
U_R=T_{z_1}e^{iR^{-1/2}\beta_1}T_{z_2}e^{iR^{-1/2}\beta_2}\cdots T_{z_{N+1}}e^{iR^{-1/2}\beta_{N+1}}\1,
\end{equation}
where the exponentials act by scalar multiplication. These multiplications are isometries on scalar $L_R$. Therefore
\begin{equation}\label{high:nested-contraction}
\int_\Omega|U_R(\omega)|^R\dd\mu(\omega)\le1.
\end{equation}
Only scalar complex contractivity is used in this construction.

\noindent\textbf{Step 2: Expansion of $U_R$ through order $1/R$.}
We give the expansion of $U_R$ through order $1/R$. For the cumulative times $Z_j=z_1+\cdots+z_j$, one has
\begin{equation}\label{high:cumulative-clock}
Z_j=2^js+\frac{i\delta\varepsilon_j2^{j+1}s}{R^{1/2}}-\frac2R\left(2s\alpha_1^2+\sum_{k=2}^j2^{k-1}s\alpha_k^2\right)+O(R^{-3/2}),
\end{equation}
which implies
\[
T_{Z_j}(\beta_j)=T_{2^js}(\beta_j)-\frac{i\delta\varepsilon_j\kappa_j}{\sqrt R}2^{j+1}sAu_{2^{j+1}s}+O(R^{-1}).
\]
At real times the product of the multiplication phases is represented by the finite path \eqref{high:path-law}. Their sum telescopes to
\[
Y=\sum_{j=1}^{N+1}\beta_j(\omega_j)=\sum_{j=1}^N\varepsilon_j\gamma_jD_j.
\]
The expansion of each phase is
\[
e^{iR^{-1/2}\beta_j}=1+\frac{i\beta_j}{R^{1/2}}-\frac{\beta_j^2}{2R}+O(R^{-3/2}).
\]
The terms of phase degree zero, one, and two in \eqref{high:nested-phase}, denoted by $\deg_0(U_R)$, $\deg_1(U_R)$, and $\deg_2(U_R)$, respectively, are
\begin{align*}
\deg_0(U_R)&=T_{z_1}\cdots T_{z_{N+1}}(\1)=\1,\\
\deg_1(U_R)&=\frac{i}{R^{1/2}}\sum_{j=1}^{N+1}T_{Z_j}(\beta_j)=\frac{i\E_\omega Y}{R^{1/2}}+\frac{\delta}{R}\sum_{j=1}^{N+1}\varepsilon_j\kappa_j2^{j+1}sAu_{2^{j+1}s}+O(R^{-3/2}),\\
\deg_2(U_R)&=-\frac1{2R}\sum_{j=1}^{N+1}T_{Z_j}(\beta_j^2)-\frac1R\sum_{1\le j<k\le N+1}T_{Z_j}(\beta_jT_{Z_k-Z_j}(\beta_k))\\
&=-\frac1{2R}\sum_{j=1}^{N+1}T_{2^js}(\beta_j^2)-\frac1R\sum_{1\le j<k\le N+1}T_{2^js}(\beta_jT_{2^ks-2^js}(\beta_k))+O(R^{-3/2})\\
&=-\frac1{2R}\sum_{j=1}^{N+1}\E_\omega[\beta_j(\omega_j)^2]-\frac1R\sum_{1\le j<k\le N+1}\E_\omega[\beta_j(\omega_j)\beta_k(\omega_k)]+O(R^{-3/2})\\
&=-\frac1{2R}\E_\omega\left(\sum_{j=1}^{N+1}\beta_j(\omega_j)\right)^2+O(R^{-3/2})=-\frac{\E_\omega Y^2}{2R}+O(R^{-3/2}).
\end{align*}
Here the functions $\beta_j$ are held fixed as the clocks vary. The real second-order correction in \eqref{high:cumulative-clock} contributes only $O(R^{-3/2})$ after multiplication by a phase of order $R^{-1/2}$; changing a clock in a term of phase degree two has the same order. Thus every clock contribution through order $R^{-1}$ is accounted for, and
\[
U_R=1+\frac{i\E_\omega Y}{\sqrt R}+\frac{4\zeta-\E_\omega Y^2}{2R}+O(R^{-3/2}),
\]
where
\[
\zeta=\frac{\delta}{2}\sum_{j=1}^{N+1}\varepsilon_j\kappa_j2^{j+1}sAu_{2^{j+1}s}
=\frac{\delta}{2}\sum_{j=1}^N\gamma_j(2^{j+1}sAu_{2^{j+1}s}-\varepsilon_j\varepsilon_{j+1}2^{j+2}sAu_{2^{j+2}s}).
\]
All these expansions hold in $L_\infty$, uniformly in $\omega$, for the fixed finite lists and fixed $\gamma_j$. This follows from the boundedness of $A$ on $L_\infty$ and of the phase multipliers.

\noindent\textbf{Step 3: An exponential estimate for fixed signs.}
Taking the squared modulus gives
\[
|U_R|^2=1+\frac{4\zeta-\E_\omega Y^2+(\E_\omega Y)^2}{R}+O(R^{-3/2}).
\]
It follows that $|U_R|^R$ converges uniformly, as $R\to\infty$, to $\exp(2\zeta-(\E_\omega Y^2-(\E_\omega Y)^2)/2)$. By \eqref{high:nested-contraction}, we obtain
\begin{align}
&\int_\Omega\exp\left[\delta\sum_{j=1}^N\gamma_j(2^{j+1}sAu_{2^{j+1}s}-\varepsilon_j\varepsilon_{j+1}2^{j+2}sAu_{2^{j+2}s})-\frac12\E_\omega\left(\sum_{j=1}^N\varepsilon_j\gamma_j(D_j-\E_\omega D_j)\right)^2\right]\dd\mu(\omega)\notag\\
&\quad=\int_\Omega\exp\left[2\zeta-\frac12(\E_\omega Y^2-(\E_\omega Y)^2)\right]\dd\mu(\omega)=\lim_{R\to\infty}\int_\Omega|U_R(\omega)|^R\dd\mu(\omega)\le1.\label{high:signed-phase}
\end{align}

\noindent\textbf{Step 4: Averaging over the signs.}
Since $\E_\varepsilon(\varepsilon_j\varepsilon_k)=\delta_{jk}$, we have
\begin{align*}
&\E_\varepsilon\sum_{j=1}^N\gamma_j(2^{j+1}sAu_{2^{j+1}s}-\varepsilon_j\varepsilon_{j+1}2^{j+2}sAu_{2^{j+2}s})=\sum_{j=1}^N\gamma_j2^{j+1}sAu_{2^{j+1}s},\\
&\E_\varepsilon\E_\omega\left[\left(\sum_{j=1}^N\varepsilon_j\gamma_j(D_j-\E_\omega D_j)\right)^2\right]
 =\sum_{j,k=1}^N\delta_{jk}\gamma_j\gamma_k\E_\omega[(D_j-\E_\omega D_j)(D_k-\E_\omega D_k)]=\sum_{j=1}^NC_j\gamma_j^2.
\end{align*}
The first identity uses the extra sign $\varepsilon_{N+1}$ also when $j=N$. Averaging \eqref{high:signed-phase} over the signs $\varepsilon$ and applying Jensen's inequality to the convex function $t\mapsto e^t$ under the expectation $\E_\varepsilon$ at each fixed $\omega$ therefore yields
\[
\int_\Omega\exp\left(\delta\sum_{j=1}^N\gamma_j2^{j+1}sAu_{2^{j+1}s}(\omega)-\frac12\sum_{j=1}^NC_j(\omega)\gamma_j^2\right)\dd\mu(\omega)\le1.
\]
Letting $\delta\uparrow1/3$ and applying Fatou's lemma proves \eqref{high:phase-bound}. This argument removes covariance between distinct times; it does not assert that the $D_j$ are conditionally orthogonal given $\omega_0$.
\end{proof}

\subsection{Gaussian integration and an exponential estimate}

We use Gaussian integration and a martingale estimate to obtain an exponential bound for $B_N(f)$ retaining the $L_2$ mass of $f$, and then deduce the restricted estimate needed for interpolation.

\begin{proposition}\label{high:exponential}
If $|f(\omega)|\le1$ almost everywhere, then
\[
\int_\Omega\left[\exp\left(\frac{B_N(f)(\omega)^2}{48}\right)-1\right]\dd\mu(\omega)\le\|f\|_{L_2(\Omega)}^2.
\]
\end{proposition}

\begin{proof}
Choose the Gaussian variance $\sigma^2=1/3$. The phase parameter has already been taken to $1/3$ in Lemma~\ref{high:phase}. Since $4/189\ge1/48$, the weighted arithmetic--geometric mean inequality gives
\begin{align*}
\int_\Omega\left[\exp\left(\frac{B_N(f)(\omega)^2}{48}\right)-1\right]\dd\mu(\omega)
&\le\frac23\int_\Omega\left[\exp\left(\frac{2B_N(f)(\omega)^2}{63}-\sum_{j=1}^N\frac{C_j(\omega)}{6}\right)-1\right]\dd\mu(\omega)\\
&\quad+\frac13\int_\Omega\left[\exp\left(\sum_{j=1}^N\frac{C_j(\omega)}{3}\right)-1\right]\dd\mu(\omega).
\end{align*}
We claim that the first integral on the right-hand side is nonpositive and that the second integral is at most $3\|f\|_{L_2(\Omega)}^2$.

For the first assertion, integrate the independent real parameters $\gamma_j$ in \eqref{high:phase-bound} against $N(0,\sigma^2)$. The phase limit has already been taken for each fixed list of parameters. Tonelli's theorem therefore justifies this integration over all $N$ time scales. For $b\in\R$ and $C\ge0$, completing the square in one dimension gives
\[
\E_\gamma\exp\left(\frac{\gamma b}{3}-\frac{C\gamma^2}{2}\right)
=\frac1{\sqrt{2\pi\sigma^2}}\int_\R\exp\left(-\frac{1+\sigma^2C}{2\sigma^2}\gamma^2+\frac{\gamma b}{3}\right)\dd\gamma
=\frac{1}{\sqrt{1+\sigma^2C}}{\exp\left(\frac{\sigma^2b^2}{18(1+\sigma^2C)}\right)},
\]
which yields
\begin{equation}\label{high:gaussian-determinant}
\int_\Omega {\prod_{j=1}^N(1+\sigma^2C_j)^{-1/2}}{\exp\left(\frac{\sigma^2}{18}\sum_{j=1}^N\frac{|2^{j+1}sAu_{2^{j+1}s}|^2}{1+\sigma^2C_j}\right)}\dd\mu(\omega)\le1.
\end{equation}
Since $|u_{2^js}(\omega_j)|\le1$, the variance identity implies
\begin{align}\label{high:trace-control}
0\le C_j(\omega)=\E_\omega D_j^2-(\E_\omega D_j)^2\le\E_\omega D_j^2\le4,\quad
\sum_{j=1}^NC_j(\omega)\le\E_\omega\sum_{j=1}^ND_j^2.
\end{align}
Thus $(1+\sigma^2C_j)^{-1}\ge3/7$ and $\log(1+\sigma^2C_j)\le\sigma^2C_j$. The integrand in \eqref{high:gaussian-determinant} is therefore at least
\[
\exp\left(\frac{2B_N(f)(\omega)^2}{63}-\sum_{j=1}^N\frac{C_j(\omega)}{6}\right),
\]
which proves the required nonpositivity.

For the second assertion, use the real-valued reverse martingale $u_{2s}(\omega_1),\ldots,u_{2^{N+1}s}(\omega_{N+1})$, and put $\mathcal F_j=\sigma(\omega_j,\ldots,\omega_{N+1})$ for $1\le j\le N+1$. Then $0\le D_j^2\le4$ and orthogonality of the increments gives
\[
\E\left[\sum_{k=1}^{j-1}D_k^2\,\middle|\,\mathcal F_j\right]
=\sum_{k=1}^{j-1}(\E[u_{2^ks}(\omega_k)^2\mid\mathcal F_j]-\E[u_{2^{k+1}s}(\omega_{k+1})^2\mid\mathcal F_j])
\le \E[u_{2s}(\omega_1)^2\mid\mathcal F_j]\le1.
\]
For $1\le k\le N$, write
\[
H_k=\sum_{1\le i_1<\cdots<i_k\le N}D_{i_1}^2\cdots D_{i_k}^2.
\]
Note that
\[
\E H_1=\sum_{j=1}^N\E D_j^2=\sum_{j=1}^N(\E u_{2^js}(\omega_j)^2-\E u_{2^{j+1}s}(\omega_{j+1})^2)\le\|u_{2s}\|_{L_2(\Omega)}^2\le\|f\|_{L_2(\Omega)}^2.
\]
Moreover, conditioning at the second index shows that for $k\ge2$,
\begin{align*}
\E H_k&=\sum_{1\le i_2<\cdots<i_k\le N}\E\left[D_{i_2}^2\cdots D_{i_k}^2\E\left[\sum_{i_1=1}^{i_2-1}D_{i_1}^2\,\middle|\,\mathcal F_{i_2}\right]\right]\\
&\le\sum_{1\le i_2<\cdots<i_k\le N}\E[D_{i_2}^2\cdots D_{i_k}^2]=\E H_{k-1}\le\cdots\le\E H_1\le\|f\|_{L_2(\Omega)}^2.
\end{align*}
Convexity of the exponential on $[0,4]$ gives $\exp(D_j^2/3)\le1+3D_j^2/4$. By \eqref{high:trace-control}, conditional Jensen's inequality and the bound on $\E H_k$ above, we obtain
\[
\int_\Omega\left[\exp\left(\sum_{j=1}^N\frac{C_j(\omega)}{3}\right)-1\right]\dd\mu(\omega)
\le\E\left[\prod_{j=1}^N\left(1+\frac34D_j^2\right)-1\right]
=\sum_{k=1}^N\left(\frac34\right)^k\E H_k\le3\|f\|_{L_2(\Omega)}^2.
\]

The two asserted integral estimates prove this proposition. The subtraction of the constant terms preserves $\|f\|_{L_2(\Omega)}^2$ on its right side.
The proof is complete.
\end{proof}

\begin{corollary}\label{high:restricted}
For $p\ge4$, every measurable $E\subset\Omega$, and every real-valued $g$ with $|g|\le1$ and $g=0$ off $E$, one has
\[
\|B_N(g)\|_{L_p(\Omega)}\le2\sqrt{3p}\,\mu(E)^{1/p}.
\]
\end{corollary}

\begin{proof}
For $p\ge4$, we have
\[
B_N(g)(\omega)^p=(12p)^{p/2}\left(\frac{B_N(g)(\omega)^2}{12p}\right)^{p/2}\le(12p)^{p/2}\left[\exp\left(\frac{B_N(g)(\omega)^2}{48}\right)-1\right].
\]
It follows from Proposition~\ref{high:exponential} that
\[
\int_\Omega B_N(g)(\omega)^p\dd\mu(\omega)\le48^{p/2}(p/4)^{p/2}\|g\|_{L_2(\Omega)}^2\le(2\sqrt{3p})^p\mu(E),
\]
which completes the proof.
\end{proof}

\subsection{Dyadic \texorpdfstring{$L_p$}{Lp} estimates for \texorpdfstring{$p\ge4$}{p >= 4}}

The purpose of this subsection is to establish $L_p$ bounds of the full dyadic operator for use in the next subsection. For every $s>0$, define $\mathcal U(f)=(2^jsAu_{2^js})_{j\in\mathbb Z}$. Note that the spectral theorem on $L_2(\Omega)$ yields
\[
\|\mathcal U(f)\|_{L_2(\Omega;\ell_2(\mathbb Z))}^2=\sum_{j\in\mathbb Z}\|2^jsAu_{2^js}\|_{L_2(\Omega)}^2
\le\sup_{\lambda>0}\left[\sum_{j\in\mathbb Z}(2^js\lambda)^2e^{-2\cdot2^js\lambda}\right]\|f\|_{L_2(\Omega)}^2\le\|f\|_{L_2(\Omega)}^2.
\]

\begin{proposition}\label{high:dyadic}
For every $p\ge4$ and every $s>0$, if $f\in L_p(\Omega)$ is real-valued, then
\[
\|\mathcal U(f)\|_{L_p(\Omega;\ell_2(\mathbb Z))}\le24\sqrt p\,\|f\|_{L_p(\Omega)}.
\]
\end{proposition}

\begin{proof}
\noindent\textbf{Step 1: Truncation.}
Suppose first that $f$ is bounded. Let
\[
d_f(t)=\mu(\{|f|>t\}),\quad g_t=\operatorname{sgn}(f)\mathbbm{1}_{\{|f|>t\}}.
\]
For $r>0$, decompose $f=f_{\le r}+f_{>r}$, where
\[
f_{\le r}=\operatorname{sgn}(f)\min\{|f|,r\}=\int_0^r g_t\dd t,\quad
f_{>r}=\operatorname{sgn}(f)(|f|-r)_+=\int_r^\infty g_t\dd t.
\]

\noindent\textbf{Step 2: A distribution-function estimate.}
For the low part, the layer sums satisfy
\[
\left\|f_{\le r}-\frac rL\sum_{\ell=1}^L g_{\ell r/L}\right\|_{L_\infty(\Omega)}\le\frac rL,\quad L\ge1.
\]
Applying Corollary~\ref{high:restricted} at exponent $3p$ to the normalized approximation error and to each layer function,
rescaling the dyadic intervals, and using monotone convergence and the triangle inequality, we obtain
\begin{align*}
\left\|\mathcal U(f_{\le r})-\frac rL\sum_{\ell=1}^L\mathcal U(g_{\ell r/L})\right\|_{L_{3p}(\Omega;\ell_2(\mathbb{Z}))}&\le12\sqrt p\,\frac rL,\\
\left\|\frac rL\sum_{\ell=1}^L\mathcal U(g_{\ell r/L})\right\|_{L_{3p}(\Omega;\ell_2(\mathbb{Z}))}&\le12\sqrt p\,\frac rL\sum_{\ell=1}^L d_f(\ell r/L)^{1/(3p)}.
\end{align*}
Since $d_f$ is nonincreasing, its bounded layer powers are Riemann integrable on $[0,r]$. Letting $L\to\infty$ yields
\[
\|\mathcal U(f_{\le r})\|_{L_{3p}(\Omega;\ell_2(\mathbb{Z}))}\le12\sqrt p\int_0^r d_f(t)^{1/(3p)}\dd t.
\]
For the high part, the $L_2$ boundedness of $\mathcal U$ gives $\|\mathcal U(f_{>r})\|_{L_2(\Omega;\ell_2(\mathbb{Z}))}\le\|(|f|-r)_+\|_{L_2(\Omega)}$.
Use these estimates for the low and high parts with $r=\lambda/(24\sqrt p)$ and split the threshold $\lambda$ into $7\lambda/8$ for $f_{\le r}$ and $\lambda/8$ for $f_{>r}$. The triangle and Chebyshev inequalities give
\begin{align}
\mu(\{\|\mathcal U(f)\|_{\ell_2(\mathbb{Z})}>\lambda\})
&\le\mu(\{\|\mathcal U(f_{\le r})\|_{\ell_2(\mathbb{Z})}>7\lambda/8\})+\mu(\{\|\mathcal U(f_{>r})\|_{\ell_2(\mathbb{Z})}>\lambda/8\})\notag\\
&\le\left(\frac8{7\lambda}\right)^{3p}\|\mathcal U(f_{\le r})\|_{L_{3p}(\Omega;\ell_2(\mathbb{Z}))}^{3p}+\left(\frac8\lambda\right)^2\|\mathcal U(f_{>r})\|_{L_2(\Omega;\ell_2(\mathbb{Z}))}^2\notag\\
&\le\left(\frac{96\sqrt p}{7\lambda}\right)^{3p}\left(\int_0^{\lambda/(24\sqrt p)}d_f(t)^{1/(3p)}\dd t\right)^{3p}\notag\\
&\quad+\frac{64}{\lambda^2}\int_\Omega\left(|f(\omega)|-\frac{\lambda}{24\sqrt p}\right)_+^2\dd\mu(\omega).\label{high:distribution}
\end{align}

\noindent\textbf{Step 3: Integration of the distribution estimate.}
We employ the weighted Hardy inequality \cite[Theorem~330]{HLP}: for nonnegative $v$, one has
\[
\int_0^\infty r^{p-\beta-1}\left(\int_0^r v(t)\dd t\right)^\beta\dd r\le\left(\frac \beta{\beta-p}\right)^\beta\int_0^\infty r^{p-1}v(r)^\beta\dd r,\quad \beta>p.
\]
Integrate \eqref{high:distribution} against $p\lambda^{p-1}\dd\lambda$. For its first term, substitute $\lambda=24\sqrt p\,r$ and apply the Hardy inequality with $\beta=3p$ and $v=d_f^{1/(3p)}$. The resulting integral is bounded by
\[
p\left(\frac{96\sqrt p}{7}\right)^{3p}(24\sqrt p)^{-2p}\int_0^\infty r^{-2p-1}\left(\int_0^r d_f(t)^{1/(3p)}\dd t\right)^{3p}\dd r
\le(6/7)^{3p}(24\sqrt p)^p\|f\|_{L_p(\Omega)}^p.
\]
For the second term, Tonelli's theorem and the same substitution give
\[
64p\int_0^\infty\lambda^{p-3}\int_\Omega\left(|f(\omega)|-\frac{\lambda}{24\sqrt p}\right)_+^2\dd\mu(\omega)\dd\lambda
=\frac{2}{9p(p-1)(p-2)}(24\sqrt p)^p\|f\|_{L_p(\Omega)}^p.
\]
Thus we have proved
\[
\|\mathcal U(f)\|_{L_p(\Omega;\ell_2(\mathbb{Z}))}^p\le\left[(6/7)^{3p}+\frac2{9p(p-1)(p-2)}\right](24\sqrt p)^p\|f\|_{L_p(\Omega)}^p
\le(24\sqrt p)^p\|f\|_{L_p(\Omega)}^p.
\]

\noindent\textbf{Step 4: Approximation.}
Bounded real-valued approximation gives the estimate on all real-valued $L_p(\Omega)$ data. Since $\mu(\Omega)=1$ and $p\ge4$, the approximation also converges in $L_2(\Omega)$, and the $L_2$ boundedness of $\mathcal U$ shows that the extension is the stated operator. This proves the proposition.
\end{proof}

\subsection{Continuous-time estimates for \texorpdfstring{$p\ge4$}{p >= 4}}

We now pass from the dyadic bound to continuous time and obtain estimates for all $p\ge4$.

\begin{proposition}\label{high:model}
For every $p\ge4$, if $f\in L_p(\Omega)$ is real-valued, then
\[
\|G_q^T(f)\|_{L_p(\Omega)}\le2^{6-3/q}\sqrt p\,\|f\|_{L_p(\Omega)},\quad 2\le q\le\infty.
\]
\end{proposition}

\begin{proof}
\noindent\textbf{Step 1: The square function bound.}
The exact logarithmic partition gives, for almost every $\omega\in\Omega$,
\[
\int_1^2\sum_{j\in\mathbb Z}|2^jsAT_{2^js}(f)(\omega)|^2\frac{\dd s}{s}=\int_0^\infty|tAT_t(f)(\omega)|^2\frac{\dd t}{t}=G_2^T(f)(\omega)^2.
\]
Since $p/2\ge1$, Minkowski's inequality in $L_{p/2}(\Omega)$ and Proposition~\ref{high:dyadic} imply
\[
\|G_2^T(f)\|_{L_p(\Omega)}\le24\sqrt{p\log2}\,\|f\|_{L_p(\Omega)}\le16\sqrt{2p}\,\|f\|_{L_p(\Omega)}.
\]

\noindent\textbf{Step 2: The derivative maximal bound.}
To treat large time exponents, keep a fixed dyadic grid. Proposition~\ref{high:dyadic} with $s=1$ gives
\begin{equation}\label{high:grid-square}
\|(2^jAu_{2^j})_{j\in \mathbb{Z}}\|_{L_p(\Omega;\ell_2(\mathbb{Z}))}\le24\sqrt p\,\|f\|_{L_p(\Omega)}.
\end{equation}
If $2^j<t<2^{j+1}$, semigroup commutation and positivity give, jointly almost everywhere in $(\omega,t)$,
\[
|tAT_t(f)(\omega)|\le2T_{t-2^j}(|2^jAu_{2^j}|)(\omega)
%\le2T_{t-2^j}(\|(2^jAu_{2^j})_{j\in \mathbb{Z}}\|_{\ell_2(\mathbb{Z})})(\omega)
\le2\sup_{t>0}T_{t}(\|(2^jAu_{2^j})_{j\in \mathbb{Z}}\|_{\ell_2(\mathbb{Z})})(\omega).
\]
For each $j$, these assertions follow first at each fixed positive time and then jointly by measurability and Fubini's theorem. The countable collection of intervals covers positive time apart from its dyadic endpoints, a set of time measure zero. Consequently Lemma~\ref{lem:heatmax} and \eqref{high:grid-square} yield
\[
\|G_\infty^T(f)\|_{L_p(\Omega)}
\le 2p'\|(2^jAu_{2^j})_{j\in\mathbb Z}\|_{L_p(\Omega;\ell_2(\mathbb Z))}\le48p'\sqrt p\,\|f\|_{L_p(\Omega)}
\le64\sqrt p\,\|f\|_{L_p(\Omega)}.
\]

\noindent\textbf{Step 3: Interpolation between the endpoint bounds.}
For $2<q<\infty$, Lemma~\ref{lem:timeinterpolation} yields
\[
\|G_q^T(f)\|_{L_p(\Omega)}\le(16\sqrt2)^{2/q}64^{1-2/q}\sqrt p\,\|f\|_{L_p(\Omega)}=2^{6-3/q}\sqrt p\,\|f\|_{L_p(\Omega)}.
\]
The same estimate at $q=2$ and $q=\infty$ follows from the endpoint bounds. This completes the proof.
\end{proof}

\section{Scalar estimates for all \texorpdfstring{$1<p<\infty$}{1 < p < infinity}}\label{sec:all}

 We combine the estimates of Sections~\ref{sec:low} and~\ref{sec:high} with the $L_2$ bound to cover the full range of $p$.

\begin{proposition}\label{prop:modelall}
For every $1<p<\infty$, if $f\in L_p(\Omega)$ is real-valued, then
\begin{equation*}
\|G_q^T(f)\|_{L_p(\Omega)}\le
\begin{cases}
\displaystyle\frac12(2^{11/2-6/q}p')^{2-p}\|f\|_{L_p(\Omega)},&1<p\le2,\\[4pt]
\displaystyle\frac12(2^{8-3/q}\sqrt p)^{2-p'}\|f\|_{L_p(\Omega)},&2\le p<\infty,
\end{cases}
\quad 2\le q\le\infty.
\end{equation*}
\end{proposition}

\begin{proof}
\noindent\textbf{Step 1: Complexification of the endpoint estimates.}
 To apply complex interpolation, we extend the estimates of Propositions~\ref{prop:low:scalar} and~\ref{high:model} to complex-valued data.
Write $g=g_1+ig_2$ with real-valued $g_1,g_2$.

For $1<p<2$, choose independent standard real Gaussian variables $\gamma_1,\gamma_2$.
If $2\le q<\infty$, Lemma~\ref{lem:gaussianmoments} and Minkowski's integral inequality with exponent $q/p\ge1$ give
\[
\E|\gamma_1|^p\,G_q^T(g)(x)^p\le \E\,G_q^T(g_1\gamma_1+g_2\gamma_2)(x)^p.
\]
For $q=\infty$, the same inequality follows from the Gaussian moment identity, since the time essential supremum of an expectation is bounded by the expectation of the time essential supremum.
Integrating in $x$, applying Proposition~\ref{prop:low:scalar}, and using Lemma~\ref{lem:gaussianmoments} once more show that
\[
\E|\gamma_1|^p\,\|G_q^T(g)\|_{L_p(\Omega)}^p\le (2^{7/2-6/q}p')^p\,\E|\gamma_1|^p\,\|g\|_{L_p(\Omega)}^p.
\]
Thus the small-$p$ estimate holds with the same constant for complex-valued $g$.

For $p\ge4$, use $|tAT_t(g)|^2=|tAT_t(g_1)|^2+|tAT_t(g_2)|^2$ and the triangle inequality in the mixed norm with exponents $p/2$ and $q/2$. Proposition~\ref{high:model} yields
\[
\|G_q^T(g)\|_{L_p(\Omega)}^2\le\|G_q^T(g_1)\|_{L_p(\Omega)}^2+\|G_q^T(g_2)\|_{L_p(\Omega)}^2\le 2(2^{6-3/q}\sqrt p)^2\|g\|_{L_p(\Omega)}^2.
\]

\noindent\textbf{Step 2: Interpolation in $p$.}
First let $2\le q<\infty$.
Apply the mixed-norm Riesz--Thorin interpolation theorem to the linear map
\[
f\longmapsto(tAT_t(f))_{t>0},\quad L_p(\Omega;\C)\longrightarrow L_p(\Omega;L_q((0,\infty),\dd t/t;\C)),
\]
in the space exponent. Lemma~\ref{lem:l2endpoint} gives an $L_2$ bound of $1/2$. If $1<p<2$, take $\theta=2-p$, so that
\[
\frac1p=\frac{1-\theta}{2}+\frac{\theta(p+1)}{2p}.
\]
The estimate at exponent $2p/(p+1)$, whose conjugate exponent is $2p'$, gives
\[
\|G_q^T(f)\|_{L_p(\Omega)}\le(1/2)^{1-\theta}(2^{7/2-6/q}\,2p')^\theta\|f\|_{L_p(\Omega)}=\frac12(2^{11/2-6/q}p')^{2-p}\|f\|_{L_p(\Omega)}.
\]
If $p>2$, take $\theta=2-p'$, so that
\[
\frac1p=\frac{1-\theta}{2}+\frac\theta{2p}.
\]
The estimate at exponent $2p\ge4$ yields
\[
\|G_q^T(f)\|_{L_p(\Omega)}\le(1/2)^{1-\theta}(2^{6-3/q}\sqrt2\sqrt{2p})^\theta\|f\|_{L_p(\Omega)}=\frac12(2^{8-3/q}\sqrt p)^{2-p'}\|f\|_{L_p(\Omega)}.
\]
The case $p=2$ follows directly from the $L_2$ bound.

Finally, for any level below $G_\infty^T(f)(x)$, the time set on which $|tAT_t(f)(x)|$ exceeds that level contains a subset of finite positive measure with respect to $\dd t/t$. Restricting the time integral to this subset gives
\[
G_\infty^T(f)(x)\le\liminf_{q\to\infty}G_q^T(f)(x).
\]
Fatou's lemma therefore gives the same estimates for $q=\infty$. This completes the proof.
\end{proof}

\section{Passage to general symmetric submarkovian semigroups}\label{sec:transfer}

The aim of this section is to prove Theorem~\ref{thm:main}.
We first transfer the real-valued estimate of Proposition~\ref{prop:modelall} from bounded conservative models to general symmetric submarkovian semigroups for $1<p<\infty$ and $2\le q\le\infty$.
Lemma~\ref{lem:transfer:finite}, inspired by \cite[(5.1)]{DKP2021} and \cite[Section 3]{Haase2016}, 
removes the conservation assumption on finite measure spaces by adjoining a single atom.
Lemma~\ref{lem:transfer:compression} uses compression to finite measure sets to obtain the estimates for $A_\varepsilon=\varepsilon^{-1}(\mathrm{Id}_{\Omega}-T_\varepsilon)$ on the original measure space.
Lemma~\ref{lem:transfer:generator} passes to the original generator as $\varepsilon\downarrow0$ and extends the estimates to all real-valued data in the corresponding $L_p$ spaces. These steps preserve the displayed constants.
In the final subsection, we pass to Hilbert-valued functions and obtain the bounds in Theorem~\ref{thm:main}.

\subsection{Removing the conservativity assumption on the operator}\label{sec5.1}

\begin{lemma}\label{lem:transfer:finite}
Let $P$ be a symmetric submarkovian kernel and its generator $A=\alpha(\mathrm{Id}_{\Omega}-P)$ for every $\alpha>0$ on a finite measure space $(\Omega,\mathcal A,\mu)$.
The estimates of Proposition~\ref{prop:modelall} hold for the semigroup $e^{-tA}$ on $(\Omega,\mathcal A,\mu)$
and every real-valued function in the corresponding $L_p$ space, with the same constants.
\end{lemma}

\begin{proof}
\noindent\textbf{Step 1: A conservative Markov extension.} Denote by $\theta(x)$ the single-step killing probability at $x$ and by $\Theta$ the total single-step mass deficit with respect to $\mu$:
\[
\theta(x)=\1-P(\1)(x),\quad \Theta=\int_\Omega\theta(x)\dd\mu(x).
\]
For $R\ge\max\{\mu(\Omega),1\}$, attach one atom $*$ of measure $R$. A function on $\Omega\cup\{*\}$ is written as $(f,z)$: it equals $f$ on $\Omega$ and has the value $z\in\R$ at $*$. Define
\begin{equation*}%
[P_R(f,z)](x)=
\begin{cases}
P(f)(x)+\theta(x)z,&x\in\Omega,\\
\dfrac1R\displaystyle\int_\Omega\theta(y)f(y)\dd\mu(y)+(1-\Theta/R)z,&x=*.
\end{cases}
\end{equation*}
More explicitly, its positive transition kernel is
\[
P_R(x,\mathrm dy)=
\begin{cases}
P(x,\mathrm dy),&x\in\Omega,\ y\in\Omega,\\[1mm]
\theta(x)\delta_*(\mathrm dy),&x\in\Omega,\ y=*,\\[1mm]
\dfrac{\theta(y)}R\,\mathrm d\mu(y),&x=*,\ y\in\Omega,\\[1mm]
\left(1-\dfrac{\Theta}{R}\right)\delta_*(\mathrm dy),&x=*,\ y=*.
\end{cases}
\]
The total transition probability from each starting point equals one since
\[
P_R(x,\Omega\cup\{*\})=
\begin{cases}
P(\1)(x)+\theta(x),&x\in\Omega,\\
\dfrac{\Theta}{R}+1-\dfrac{\Theta}{R},&x=*.
\end{cases}
\]
Symmetry on the $\Omega$ block follows from that of $P$, and the two off-diagonal terms agree because
\[
R[P_R(f,0)](*)z
=\int_\Omega\theta(x)f(x)z\dd\mu(x)=\int_\Omega f(x)[P_R(0,z)](x)\dd\mu(x).
\]
Thus $P_R$ is a symmetric Markov kernel on $(\Omega\cup\{*\},\mu+R\delta_*)$. Normalizing $\mu+R\delta_*$ by its total mass $\mu(\Omega)+R$ gives a probability measure. Both sides of each estimate scale by $(\mu(\Omega)+R)^{-1/p}$, so the estimates also hold for the unnormalized measure, with unchanged constants. For fixed $\alpha>0$, the bounded-model estimates apply to $\alpha(\mathrm{Id}_{\Omega\cup\{*\}}-P_R)$ and the datum $(f,0)$, whose norm is exactly
\[
\|(f,0)\|_{L_p(\Omega\cup\{*\},\mu+R\delta_*)}=\|f\|_{L_p(\Omega,\mu)}.
\]

\noindent\textbf{Step 2: $L_p$ convergence of the restricted powers.} For bounded real-valued $f$, write $P_R^n(f,0)=(f_{R,n},z_{R,n})$. The two components obey the recursions
\[
\begin{aligned}
f_{R,0}&=f,&z_{R,0}&=0,\\
f_{R,n+1}&=P(f_{R,n})+\theta\,z_{R,n},& z_{R,n+1}&=\frac1R\int_\Omega\theta(x)f_{R,n}(x)\dd\mu(x)+(1-\Theta/R)z_{R,n}.
\end{aligned}
\]
Contractivity gives $\|f_{R,n}\|_{L_\infty(\Omega)}\le\|f\|_{L_\infty(\Omega)}$. Iteration of the first recursion, and then of the second, gives
\[
f_{R,n}-P^n(f)=\sum_{\ell=0}^{n-1}P^{n-1-\ell}(\theta\,z_{R,\ell}),\quad z_{R,n}=\frac1R\sum_{\ell=0}^{n-1}\left(1-\frac{\Theta}{R}\right)^{n-1-\ell}\int_\Omega\theta(x)f_{R,\ell}(x)\,\mathrm d\mu(x),
\]
and therefore\footnote{\ Here we first estimate $|z_{R,n}|$ and then use this bound to estimate $\|f_{R,n}-P^n(f)\|_{L_p(\Omega)}$.}
\[
\|f_{R,n}-P^n(f)\|_{L_p(\Omega)}\le\frac{n(n-1)\mu(\Omega)^{1+1/p}}{2R}\|f\|_{L_\infty(\Omega)},\quad |z_{R,n}|\le\frac{n\mu(\Omega)}R\|f\|_{L_\infty(\Omega)}.
\]
For arbitrary real-valued $f\in L_p(\Omega)$, choose bounded real-valued $g\in L_p(\Omega)$ such that $\|f-g\|_{L_p(\Omega)}$ is arbitrarily small. The contractions $P_R^n$ and $P^n$ give
\[
\|[P_R^n(f,0)]|_\Omega-P^n(f)\|_{L_p(\Omega)} \le2\|f-g\|_{L_p(\Omega)}+\frac{n(n-1)\mu(\Omega)^{1+1/p}}{2R}\|g\|_{L_\infty(\Omega)}.
\]
Letting $R\to\infty$ proves convergence for every fixed $n$ and $p$.

\noindent\textbf{Step 3: Uniform $L_p$ convergence of the first-order trajectories on bounded time intervals.} Write $A_R=\alpha(\mathrm{Id}_{\Omega\cup\{*\}}-P_R)$. The Poisson expansion is
\begin{equation}\label{transfer:poisson}
A_Re^{-tA_R}(f,0)=\alpha e^{-\alpha t}\sum_{n=0}^\infty\frac{(\alpha t)^n}{n!}(\mathrm{Id}_{\Omega\cup\{*\}}-P_R)P_R^n(f,0).
\end{equation}
For $0\le t\le\tau$, the norm of the part with $n>N$, after multiplication by $t$, is at most
\[
2\alpha\tau\sum_{n>N}\frac{(\alpha\tau)^n}{n!}\|f\|_{L_p(\Omega)}\longrightarrow0,\quad N\to\infty,
\]
uniformly in $R$. For the finitely many remaining terms, use the fixed-power convergence of $P_R^n(f,0)$ and $P_R^{n+1}(f,0)$ after restriction to $\Omega$. This gives the  trajectory limit, with the operators acting on $(f,0)$ on the enlarged space:
\[
\lim_{R\to\infty}\sup_{0\le t\le\tau} \left(\int_\Omega |[tA_Re^{-tA_R}(f,0)](x)-tAe^{-tA}f(x)|^p\dd\mu(x)\right)^{1/p}=0.
\]
The same Poisson expansion and tail estimate hold for $P$ acting on $f$ in $L_p(\Omega)$.

\noindent\textbf{Step 4: Pointwise control of $g$-functions for $A$ by those for $A_R$.} To pass the square-function and maximal estimates, note that for every $0<a<b<\infty$,
\begin{align*}
&\int_a^b\int_\Omega|[tA_Re^{-tA_R}(f,0)](x)-tAe^{-tA}(f)(x)|^p\dd\mu(x)\frac{\dd t}{t}\\
&\quad\le\log(b/a)\sup_{a\le t\le b}\|[tA_Re^{-tA_R}(f,0)]|_\Omega-tAe^{-tA}(f)\|_{L_p(\Omega)}^p\longrightarrow0.
\end{align*}
Choose a subsequence along which $[tA_Re^{-tA_R}(f,0)](x)$ converges to $tAe^{-tA}(f)(x)$ almost everywhere on $\Omega\times(0,\infty)$, using the windows $[1/m,m]$, $m=2,3,\ldots$. For $2\le q<\infty$, two applications of Fatou's lemma give
\[
\begin{aligned}
\int_\Omega\left(\int_a^b|tAe^{-tA}(f)(x)|^q\frac{\dd t}{t}\right)^{p/q}\dd\mu(x) \le\liminf_{R\to\infty}\int_\Omega\left(\int_a^b|[tA_Re^{-tA_R}(f,0)](x)|^q\frac{\dd t}{t}\right)^{p/q}\dd\mu(x),
\end{aligned}
\]
where here and below $R$ runs through that subsequence. At $q=\infty$, the required lower semicontinuity is
\begin{equation}\label{transfer:esssup}
\esssup_{a<t<b}|tAe^{-tA}(f)(x)| \le\liminf_{R\to\infty}\esssup_{a<t<b}|[tA_Re^{-tA_R}(f,0)](x)|,
\end{equation}
for almost every $x\in\Omega$. Indeed, for almost every $x$ the convergence holds for almost every $t$; take a subsequence realizing the lower limit on the right and then the essential supremum in $t$. Spatial Fatou's lemma applies to \eqref{transfer:esssup}. Letting $a\downarrow0$ and $b\uparrow\infty$ passes the full bounds. The integral over the added atom is nonnegative, so restricting the estimates to $\Omega$ preserves their constants.

\noindent\textbf{Step 5: The desired $L_p$ estimates.}
Apply Proposition~\ref{prop:modelall} to $A_R$ and $(f,0)$ on the enlarged conservative space. The restriction to $\Omega$, the norm identity in Step 1, and the limiting argument in Step 4 give the desired bounds for the semigroup $e^{-tA}$, with unchanged constants.
The proof is complete.
\end{proof}

\begin{remark}
Lemma~\ref{lem:transfer:finite} also admits a completely different direct proof, which neither adds an extra atom nor uses approximation. Equip $\widehat\Omega=\Omega\times\{-1,1\}$ with the measure
\[
\widehat\mu =\mu\otimes\frac12(\delta_{-1}+\delta_1).
\]
Set $\theta=\1-P(\1)$, let $M_\theta$ denote multiplication by $\theta$, and define
\[
\widehat P=
\begin{pmatrix}
P+\dfrac12M_\theta&\dfrac12M_\theta\\
\dfrac12M_\theta&P+\dfrac12M_\theta
\end{pmatrix}.
\]
Then $\widehat P$ is a measurable symmetric Markov kernel. Indeed, its entries are positive, its row sums preserve $\1$, and symmetry follows from that of $P$ and $M_\theta$.

The map $J(f)(x,\pm1)=\pm f(x)$ is an isometry from $L_p(\Omega)$ into $L_p(\widehat\Omega)$. The added multiplication terms cancel on its range, giving $\widehat PJ=JP$. Consequently, with this notation
\[
A=\alpha(\mathrm{Id}_{\Omega}-P),\quad \widehat A =\alpha(\mathrm{Id}_{\widehat\Omega}-\widehat P), \quad T_t=e^{-tA},\quad \widehat T_t=e^{-t\widehat A},
\]
we have
\[
\widehat T_tJ=JT_t,\quad t\widehat A\widehat T_tJ=JtAT_t.
\]
Thus, for every $2\le q\le\infty$,
\[
\|G_q^{\widehat T}(Jf)\|_{L_p(\widehat\Omega)} = \|G_q^T(f)\|_{L_p(\Omega)}.
\]
Applying Proposition~\ref{prop:modelall} to $Jf$ therefore proves the lemma with exactly the same constants. Here the finite measure $\widehat\mu$ may first be normalized to a probability measure; undoing this normalization changes both sides by the same factor.
\end{remark}

\subsection{Removing the finiteness assumption on the measure}

\begin{lemma}\label{lem:transfer:compression}
For every $\varepsilon>0$, let $A_\varepsilon=(\mathrm{Id}_{\Omega}-T_\varepsilon)/\varepsilon$.
The estimates of Proposition~\ref{prop:modelall} hold for the semigroup $e^{-tA_\varepsilon}$ on $(\Omega,\mathcal A,\mu)$
and every real-valued function in the corresponding $L_p$ space, with the same constants.
\end{lemma}

\begin{proof}
Choose measurable sets $\Omega_m\uparrow\Omega$ with $\mu(\Omega_m)<\infty$, and define $P_m(g)=\mathbbm{1}_{\Omega_m}T_\varepsilon(\mathbbm{1}_{\Omega_m}g)$ on $\Omega$. Let $A_{\varepsilon,m}=\varepsilon^{-1}(\mathrm{Id}_{\Omega}-P_m)$ on $L_p(\Omega)$. Since $P_m$ vanishes on functions supported in $\Omega\setminus\Omega_m$, we have the block decomposition
\[
A_{\varepsilon,m}
=\underbrace{\varepsilon^{-1}(\mathrm{Id}_{\Omega_m}-P_m|_{L_p(\Omega_m)})}_{\text{on }L_p(\Omega_m)}
\oplus \underbrace{\varepsilon^{-1}\mathrm{Id}_{\Omega\setminus\Omega_m}}_{\text{on }L_p(\Omega\setminus\Omega_m)}.
\]
The restriction of $P_m$ to $\Omega_m$ is a symmetric submarkovian kernel, so Lemma~\ref{lem:transfer:finite} applies to the first block.
For a function supported in $\Omega_m$, the semigroup and its derivatives therefore coincide with the finite-measure construction on $\Omega_m$ and vanish on the other block.

For every $g\in L_p(\Omega)$ with $p<\infty$, we have
\begin{align*}
\|P_m(g)-T_\varepsilon(g)\|_{L_p(\Omega)}
&\le \|T_\varepsilon(\mathbbm{1}_{\Omega\setminus\Omega_m}g)\|_{L_p(\Omega)}
   +\|\mathbbm{1}_{\Omega\setminus\Omega_m}T_\varepsilon(g)\|_{L_p(\Omega)}\\
&\le \|\mathbbm{1}_{\Omega\setminus\Omega_m}g\|_{L_p(\Omega)}
   +\|\mathbbm{1}_{\Omega\setminus\Omega_m}T_\varepsilon(g)\|_{L_p(\Omega)}
   \longrightarrow0.
\end{align*}
For each fixed integer $j\ge1$, the telescoping identity
\[
(P_m^j-T_\varepsilon^j)g =\sum_{\ell=0}^{j-1}P_m^{j-1-\ell}(P_m-T_\varepsilon)T_\varepsilon^\ell(g)
\]
and contractivity yield
\[
\|(P_m^j-T_\varepsilon^j)g\|_{L_p(\Omega)} \le\sum_{\ell=0}^{j-1}\|(P_m-T_\varepsilon)T_\varepsilon^\ell(g)\|_{L_p(\Omega)} \longrightarrow0,\quad m\to \infty.
\]
Apply \eqref{transfer:poisson} and its uniform tail estimate to $P_m$ and $T_\varepsilon$. With $A_\varepsilon=(\mathrm{Id}_{\Omega}-T_\varepsilon)/\varepsilon$, this gives for every $f\in L_p(\Omega)$ and $0<\tau<\infty$,
\[
\lim_{m\to\infty}\sup_{0\le t\le\tau}\|tA_{\varepsilon,m}e^{-tA_{\varepsilon,m}}f-tA_\varepsilon e^{-tA_\varepsilon}f\|_{L_p(\Omega)}=0.
\]
For arbitrary real-valued $f\in L_p(\Omega)$, use $\mathbbm{1}_{\Omega_m}f$ in the finite-measure estimates. Since $\|\mathbbm{1}_{\Omega_m}f\|_{L_p(\Omega)}\le\|f\|_{L_p(\Omega)}$, $\|\mathbbm{1}_{\Omega_m}f-f\|_{L_p(\Omega)}\to0$ and
\[
\sup_{0\le t\le\tau}\|tA_{\varepsilon,m}e^{-tA_{\varepsilon,m}}(\mathbbm{1}_{\Omega_m}f-f)\|_{L_p(\Omega)}\le(2\tau/\varepsilon)\|\mathbbm{1}_{\Omega_m}f-f\|_{L_p(\Omega)},
\]
the same trajectory convergence holds for these varying data. The preceding arguments using Fatou's lemma pass all the required estimates to $A_\varepsilon$, with unchanged constants.
\end{proof}

\subsection{Removing the boundedness assumption on the generator}

\begin{lemma}\label{lem:transfer:generator}
The estimates of Proposition~\ref{prop:modelall} hold for the original semigroup $T_t$ on $(\Omega,\mathcal A,\mu)$
and every real-valued function in the corresponding $L_p$ space, with the same constants.
\end{lemma}

\begin{proof}
Fix $1<p<\infty$ and $2\le q\le\infty$. By Lemma~\ref{lem:transfer:compression}, the required estimates hold for $e^{-tA_\varepsilon}$ uniformly in $\varepsilon>0$.

First let $f\in L_p(\Omega)\cap L_2(\Omega)$ be real-valued. On $L_2(\Omega)$, spectral calculus gives
\[
A_\varepsilon=a_\varepsilon(A), \quad a_\varepsilon(\lambda) =\frac{1-e^{-\varepsilon\lambda}}{\varepsilon}.
\]
For every $\lambda\ge0$, $ta_\varepsilon(\lambda)e^{-ta_\varepsilon(\lambda)}$ converges to $t\lambda e^{-t\lambda}$ locally uniformly for $t>0$, and both multipliers are bounded by $1/e$. The spectral theorem and dominated convergence yield
\[
\lim_{\varepsilon\downarrow0} \sup_{a\le t\le b} \|tA_\varepsilon e^{-tA_\varepsilon}f-tAT_t(f) \|_{L_2(\Omega)}=0, \quad 0<a<b<\infty.
\]
The diagonal subsequence and Fatou arguments in Step~4 of the proof of Lemma~\ref{lem:transfer:finite} now pass the estimates to $T_t$, with unchanged constants; the case $q=\infty$ follows using \eqref{transfer:esssup}.

Finally, note that $L_p(\Omega)\cap L_2(\Omega)$ is dense in $L_p(\Omega)$. The local derivative bounds \eqref{eq:localderivatives} and the same Fatou argument extend the estimate to every real-valued $f\in L_p(\Omega)$. The proof is complete.
\end{proof}

\subsection{From scalar-valued to Hilbert-valued functions}\label{sec:hilbert}

\begin{proof}[Proof of Theorem~\ref{thm:main}]
Fix $1<p<\infty$ and $2\le q\le\infty$. Lemma~\ref{lem:transfer:generator}
gives the real-valued estimate of Proposition~\ref{prop:modelall}
for the original semigroup. We now pass to Hilbert-valued functions.

\noindent\textbf{Step 1: Finite-dimensional real Hilbert spaces.}
First suppose that $\mathcal{H}$ is a finite-dimensional real Hilbert space with an orthonormal basis $(e_k)_{k=1}^n$, and let $f=\sum_{k=1}^n f_ke_k$ belong to $L_p(\Omega;\mathcal{H})$. Choose independent standard real Gaussian variables $\gamma_1,\ldots,\gamma_n$.

The scalar estimate of Proposition~\ref{prop:modelall},
transferred to the original setting by Lemma~\ref{lem:transfer:generator},
yields a bounded real linear operator
\begin{align*}
S:L_p(\Omega)&\longrightarrow L_p(\Omega;L_q((0,\infty),\dd t/t;\R)),\\
g&\longmapsto Sg=(tAT_t(g))_{t>0},
\end{align*}
whose norm satisfies the bounds in Proposition~\ref{prop:modelall}.
At $q=\infty$, the target is understood as the space of jointly measurable functions with finite mixed norm, using the essential supremum in time. Define
\[
S^*:L_{p'}(\Omega;L_{q'}((0,\infty),\dd t/t;\R))\longrightarrow L_{p'}(\Omega)
\]
by the integral pairing
\[
\int_\Omega g(x)S^*h(x)\dd\mu(x)=\int_\Omega\int_0^\infty Sg(x,t)h(x,t)\frac{\dd t}{t}\dd\mu(x).
\]
For $2\le q<\infty$, this is the Banach space adjoint of $S$. At $q=\infty$, we take $q'=1$, and $S^*$ is the restriction of the Banach space adjoint to the integral test space above. In both cases the test space is norming for the mixed norm, so $\|S^*\|=\|S\|$.

Let $h_1,\ldots,h_n\in L_{p'}(\Omega;L_{q'}((0,\infty),\dd t/t;\R))$. Applying Lemma~\ref{lem:gaussianmoments} to $a_k=S^*h_k$ and then using the bound on $S^*$, we obtain
\begin{align*}
\left\|\left(\sum_{k=1}^n|S^*h_k|^2\right)^{1/2}\right\|_{L_{p'}(\Omega)}
&=
\frac{1}{(\E|\gamma_1|^{p'})^{1/p'}}
\left(\E\left\|S^*\left(\sum_{k=1}^n h_k\gamma_k\right)\right\|_{L_{p'}(\Omega)}^{p'}\right)^{1/p'}\\
&\le
\frac{\|S^*\|}{(\E|\gamma_1|^{p'})^{1/p'}}
\left(\E\left\|\sum_{k=1}^n h_k\gamma_k\right\|_{L_{p'}(\Omega;L_{q'}((0,\infty),\dd t/t;\R))}^{p'}\right)^{1/p'}\\
&\le
\|S\|
\frac{(\E|\gamma_1|^{\max\{p',q'\}})^{1/\max\{p',q'\}}}
     {(\E|\gamma_1|^{p'})^{1/p'}}
\left\|\left(\sum_{k=1}^n|h_k|^2\right)^{1/2}\right\|_{L_{p'}(\Omega;L_{q'}((0,\infty),\dd t/t;\R))}.
\end{align*}
For the last inequality, use Tonelli's theorem to interchange the spatial integral and the Gaussian expectation, and then use monotonicity of Gaussian moments to pass from $p'$ to $\max\{p',q'\}$.
Since $\max\{p',q'\}\ge q'$, Minkowski's inequality moves the Gaussian norm inside the time norm, and Lemma~\ref{lem:gaussianmoments} gives the pointwise square sum.
The Gaussian moment ratio is $1$ if $p'\ge2$, since $q'\le2$. If $p'<2$, its numerator is at most $(\E|\gamma_1|^2)^{1/2}=1$. Gaussian integration by parts gives $\E|\gamma_1|^{p'+2}=(p'+1)\E|\gamma_1|^{p'}$. By H\"{o}lder's inequality and $p'+1\le2^{p'}$,
\[
1=\E|\gamma_1|^2\le(\E|\gamma_1|^{p'})^{p'/2}(\E|\gamma_1|^{p'+2})^{1-p'/2}
\le(\E|\gamma_1|^{p'})\,2^{p'(2-p')/2}.
\]
Thus the moment ratio is at most $2^{(2-p')/2}$ if $p'<2$,
and hence at most $2^{\max\{2-p',0\}/2}$ for all $1<p<\infty$.

For $\mathcal{H}$-valued $f$, use the same notation $S(f)(x,t)=tAT_t(f)(x)$,
and let $h=\sum_{k=1}^n h_ke_k$ range over the unit ball of $L_{p'}(\Omega;L_{q'}((0,\infty),\dd t/t;\mathcal{H}))$.
This test space is norming for the mixed norm. For $q=\infty$, the same conclusion follows by choosing measurable maximizers over finite subsets of a countable dense set in the time $L_1$ unit ball and then applying $L_p$ duality in $\Omega$.
By the definition of $S^*$, Cauchy--Schwarz's inequality in $\mathcal{H}$, H\"{o}lder's inequality in $\Omega$,
and the preceding estimate, we obtain
\begin{align*}
\|G_q^T(f)\|_{L_p(\Omega)}
&=\sup_h\left|\int_\Omega\int_0^\infty\langle S(f)(x,t),h(x,t)\rangle_{\mathcal{H}}\frac{\dd t}{t}\dd\mu(x)\right|\\
&\le\sup_h\left\|\left(\sum_{k=1}^n|S^*h_k|^2\right)^{1/2}\right\|_{L_{p'}(\Omega)} \|f\|_{L_p(\Omega;\mathcal{H})} \\
&\le 2^{\max\{2-p',0\}/2}\|S\|\|f\|_{L_p(\Omega;\mathcal{H})} \\
&\le
 \begin{cases}
 \displaystyle\frac12{(2^{11/2-6/q}p')}^{2-p}\|f\|_{L_p(\Omega;\mathcal{H})},
       &1<p\le2,\\[4pt]
 \displaystyle\frac12{(2^{17/2-3/q}\sqrt p)}^{2-p'}\|f\|_{L_p(\Omega;\mathcal{H})},
      &2\le p<\infty,
 \end{cases}
\end{align*}
which proves Theorem~\ref{thm:main} for finite-dimensional real Hilbert spaces.

\noindent\textbf{Step 2: Passage to general Hilbert spaces.}
Let $\mathcal{H}$ be a real Hilbert space and let $f\in L_p(\Omega;\mathcal{H})$. By strong measurability, $f$ takes values almost everywhere in a separable closed subspace of $\mathcal{H}$. Let $\pi_m$ be increasing finite-rank orthogonal projections onto subspaces of this closed subspace, converging strongly to the identity on it.
The kernel action of $T_t$ preserves functions taking values in this closed subspace and commutes with $\pi_m$. By differentiation, $tAT_t(f)$ also takes values in this subspace, and $tAT_t(\pi_mf)=\pi_m(tAT_t(f))$.
Thus $\|tAT_t(\pi_mf)(x)\|_{\mathcal{H}}\uparrow\|tAT_t(f)(x)\|_{\mathcal{H}}$ for almost every $(x,t)$. For $2\le q<\infty$, monotone convergence in time gives $G_q^T(\pi_mf)\uparrow G_q^T(f)$ almost everywhere. The same conclusion holds for $q=\infty$, since the essential supremum in time commutes with the countable supremum over $m$.
Apply the finite-dimensional estimate of Step 1 to $\pi_mf$ and use $\|\pi_mf\|_{L_p(\Omega;\mathcal{H})}\le\|f\|_{L_p(\Omega;\mathcal{H})}$. Monotone convergence in $\Omega$ gives the same estimate for $f$. A complex Hilbert space may be regarded as a real Hilbert space with the same norm, and the real-time kernel action is unchanged.

Finally, note that $AT_t(F(f))=0$ and therefore $G_q^T(f-F(f))=G_q^T(f)$.
Applying the estimate just proved to $f-F(f)$ gives \eqref{eq:main}. The proof is complete.
\end{proof}

\section{Sharpness for the Euclidean heat semigroup}\label{sec:sharp}

Let $H^{(d)}_t=e^{t\Delta}$ be the heat semigroup on $\R^d$.
Its fixed point projection is zero on
$L_p(\R^d;\mathcal{H})$ for $1<p<\infty$. Indeed, heat kernel estimates give
strong convergence to zero at large times for bounded compactly
supported data, and contractivity and density give this convergence
for every $L_p$ function. Thus the constants in \eqref{eq:best} are the
operator norms with denominator $\|f\|_{L_p(\R^d;\mathcal{H})}$ in this section.
We first obtain a lower bound by
constructing compactly supported approximate eigenfunctions in growing
dimensions. The approximation will be used on a finite time interval
before its length tends to infinity.

\begin{proposition}\label{sharp:shell}
For $1<p<\infty$, one has
\begin{equation}\label{sharp:eq-shell}
 \sup_{d\ge1}\mathrm L_{c,q,p}^{H^{(d)}}(\C)
 \ge
 \begin{cases}
 \displaystyle\frac{\Gamma(q)^{1/q}}{q}
                  \frac{p}{2\sqrt{p-1}},&2\le q<\infty,\\[4pt]
 \displaystyle\frac{p}{2e\sqrt{p-1}},&q=\infty.
 \end{cases}
\end{equation}
The same lower bounds hold for any nonzero complex Hilbert space and
for any real Hilbert space of dimension at least two. For real-valued
functions, and hence for any nonzero real Hilbert space, the right-hand
sides in \eqref{sharp:eq-shell} may be divided by $2$.
\end{proposition}

\begin{proof}
Fix $p$ and put
\[
 s=-\frac1p+i\frac{\sqrt{p-1}}p,\quad
 \lambda=-s^2-s=\frac{2(p-1)}{p^2}
                  +i\frac{(2-p)\sqrt{p-1}}{p^2}.
\]
Let
$0\ne\chi\in C_c^\infty((-2,2))$ and, for positive integers $n$, set
\[
 f_n(x)=e^{s(|x|-n)}
              \chi((|x|-n)/n^{1/3}),\quad x\in\R^{n+1}.
\]
For all sufficiently large $n$, the support of $f_n$ is contained in
$n-2n^{1/3}<|x|<n+2n^{1/3}$ and is separated from the origin. In particular,
$f_n\in C_c^\infty(\R^{n+1};\C)$ and belongs to $\mathrm{Dom}(-\Delta)$.

In the usual radial coordinate $r=|x|$, the measure is
$|\mathbb S^n|r^n\dd r$, where $|\mathbb S^n|$ is the surface area
of the unit sphere in $\R^{n+1}$. Taylor's formula gives, for large $n$,
\[
 \sup_{-2\le u\le2} |n\log(1+n^{-2/3}u)-n^{1/3}u| \le8n^{-1/3}.
\]
Since $p\operatorname{Re}s=-1$, the exponential factors cancel in the
$p$th power of the norm, and
\begin{equation}\label{sharp:eq-volume}
 \|f_n\|_{L_p(\R^{n+1})}^p  =|\mathbb S^n|n^{n+1/3} \int_{-2}^2|\chi(u)|^p    e^{n\log(1+n^{-2/3}u)-n^{1/3}u}\dd u.
\end{equation}
The last integral tends to $\|\chi\|_{L_p(-2,2)}^p>0$.

On radial functions in $\R^{n+1}$, the Laplacian acts as
\[
 \Delta=\partial_r^2+\frac nr\partial_r,\quad r>0.
\]
With $L_\lambda=-\Delta-\lambda$, we have
\begin{equation}\label{sharp:eq-residual}
 L_\lambda(f_n)(x)=e^{s(r-n)}\left[
    \left(1-\frac nr\right)s\chi\left(\frac{r-n}{n^{1/3}}\right)
    -\frac{2s+n/r}{n^{1/3}}\chi'\left(\frac{r-n}{n^{1/3}}\right)
    -\frac1{n^{2/3}}\chi''\left(\frac{r-n}{n^{1/3}}\right)
 \right].
\end{equation}
On the support, $n/r\le2$ and $|1-n/r|\le4n^{-2/3}$ for large $n$.
Applying the same measure computation as in \eqref{sharp:eq-volume}
to the three terms in \eqref{sharp:eq-residual} gives
\begin{align*}
\varepsilon_n
 &=\frac{\|L_\lambda(f_n)\|_{L_p(\R^{n+1})}} {\|f_n\|_{L_p(\R^{n+1})}} \le 4|s|n^{-2/3}
 +(2|s|+2)n^{-1/3} \frac{ \left\|e^{s(r-n)} \chi'\left(\frac{r-n}{n^{1/3}}\right) \right\|_{L_p(\R^{n+1})}} {\|f_n\|_{L_p(\R^{n+1})}}\\
 &\quad+n^{-2/3} \frac{ \left\|e^{s(r-n)} \chi''\left(\frac{r-n}{n^{1/3}}\right) \right\|_{L_p(\R^{n+1})}} {\|f_n\|_{L_p(\R^{n+1})}}\\
 &\le 4|s|n^{-2/3} +e^{16n^{-1/3}/p} \left[ (2|s|+2)n^{-1/3} \frac{\|\chi'\|_{L_p(-2,2)}}{\|\chi\|_{L_p(-2,2)}} +n^{-2/3}
   \frac{\|\chi''\|_{L_p(-2,2)}}{\|\chi\|_{L_p(-2,2)}} \right]\\
 &\le C_{p,\chi}\left(n^{-1/3}+n^{-2/3}\right)
  \longrightarrow0,\quad n\to\infty.
\end{align*}
No uniformity in $p$ is required in this approximation.

For the rest of the proof, write $H_t=H^{(n+1)}_t$. Duhamel's formula
in $L_p$ yields
\[
 H_t(f_n)-e^{-\lambda t}f_n
 =-\int_0^t e^{-\lambda(t-u)}H_u(L_\lambda(f_n))\dd u.
\]
Contractivity and $\operatorname{Re}\lambda>0$ imply
$\|H_t(f_n)-e^{-\lambda t}f_n\|_{L_p(\R^{n+1})}\le t\|L_\lambda(f_n)\|_{L_p(\R^{n+1})}$,
and hence by the fact that $-\Delta$ commutes with $H_t$ on $\mathrm{Dom}(-\Delta)$
\begin{equation}\label{sharp:eq-derivative}
 \|(-\Delta)H_t(f_n)-\lambda e^{-\lambda t}f_n\|_{L_p(\R^{n+1})}
 \le(1+|\lambda|t)\|L_\lambda(f_n)\|_{L_p(\R^{n+1})}.
\end{equation}
Thus the derivative itself is controlled, without differentiating an
estimate for the semigroup error.

Choose a norm-one complex linear functional $\ell_n$ on $L_p(\R^{n+1};\C)$ such that $\ell_n(f_n)=\|f_n\|_{L_p(\R^{n+1})}$. Fix $0<T<\infty$. The representation of $\ell_n$ by an $L_{p'}(\R^{n+1};\C)$ function, Minkowski's integral inequality in time, and H\"{o}lder's inequality in space give
\[
\bigl\|\ell_n(t(-\Delta)H_t(f_n))\bigr\|_{L_q((0,T),\dd t/t;\C)}
\le \|G_q^{H^{(n+1)}}(f_n)\|_{L_p(\R^{n+1})},
\quad 2\le q\le\infty,
\]
where the time norm at $q=\infty$ is the essential supremum. It follows from  \eqref{sharp:eq-derivative} that
\[
\left|\ell_n(t(-\Delta)H_t(f_n))
-t\lambda e^{-\lambda t}\|f_n\|_{L_p(\R^{n+1})}\right|
\le t(1+|\lambda|t)\|L_\lambda(f_n)\|_{L_p(\R^{n+1})}.
\]
Taking the time norm over $(0,T)$ and using $\|t\|_{L_q((0,T),\dd t/t)}\le T$, we obtain
\[
\frac{\|G_q^{H^{(n+1)}}(f_n)\|_{L_p(\R^{n+1})}}
     {\|f_n\|_{L_p(\R^{n+1})}}
\ge
\|t\lambda e^{-\lambda t}\|_{L_q((0,T),\dd t/t;\C)}
-T(1+|\lambda|T)\varepsilon_n.
\]
Letting first $n\to\infty$ and then $T\to\infty$ gives
\[
\sup_{d\ge1}\mathrm L_{c,q,p}^{H^{(d)}}(\C)
\ge \|t\lambda e^{-\lambda t}\|_{L_q((0,\infty),\dd t/t;\C)}
=
\begin{cases}
 \displaystyle\frac{\Gamma(q)^{1/q}}{q}
                  \frac{p}{2\sqrt{p-1}},&2\le q<\infty,\\[4pt]
 \displaystyle\frac{p}{2e\sqrt{p-1}},&q=\infty,
 \end{cases}
\]
which proves \eqref{sharp:eq-shell}.

Finally, the identification $\C\cong\R^2$ is isometric and intertwines
the componentwise heat semigroups. The complex scalar norms above
are therefore exactly those for the real two-dimensional Hilbert
space. Isometric embeddings give the asserted lower bounds for larger
real Hilbert spaces and for nonzero complex Hilbert spaces. For the
real scalar case, write $f=u+iv$. For every $2\le q\le\infty$,
\[
 G_q^{H^{(d)}}(f)
 \le G_q^{H^{(d)}}(u)+G_q^{H^{(d)}}(v),\quad
 \|u\|_{L_p(\R^d)},\ \|v\|_{L_p(\R^d)}\le\|f\|_{L_p(\R^d)}.
\]
It follows that
$\mathrm L_{c,q,p}^{H^{(d)}}(\C)\le2\mathrm L_{c,q,p}^{H^{(d)}}(\R)$.
Embedding the real scalars into any nonzero real Hilbert space completes
the proof.
\end{proof}

\begin{proposition}\label{sharp:tail}
For every $1<p\le2$ and $2\le q\le\infty$,
\[
 \mathrm L_{c,q,p}^{H^{(1)}}(\R)\ge\frac{p'}{512}.
\]
The estimate is witnessed by one fixed nonnegative function in
$C_c^\infty((0,1))$, independent of both $p$ and $q$.
\end{proposition}

\begin{proof}
Choose $\eta\in C_c^\infty((0,1))$ with $0\le\eta\le1$ and
$\int_0^1\eta(y)\dd y\ge1/2$. In particular, $\|\eta\|_{L_p(\R)}\le1$.
The one-dimensional heat kernel gives
\[
 t(-\Delta)H_t^{(1)}\eta(x)
 =\int_0^1\left(\frac12-\frac{(x-y)^2}{4t}\right)
     (4\pi t)^{-1/2}e^{-(x-y)^2/(4t)}\eta(y)\dd y.
\]
For $x\ge2$ and $x^2\le t\le2x^2$, we have
$(x-y)^2/(4t)\le1/4$. All integrands are nonnegative, and
\[
 t(-\Delta)H_t^{(1)}\eta(x)\ge\frac{e^{-1/4}}{16\sqrt{2\pi}\,x}.
\]
The interval $[x^2,2x^2]$ has $\dd t/t$ measure $\log2$. Thus
\[
 G_q^{H^{(1)}}(\eta)(x)
 \ge\frac{e^{-1/4}(\log2)^{1/q}}{16\sqrt{2\pi}\,x}
 \ge\frac{e^{-1/4}\sqrt{\log2}}{16\sqrt{2\pi}\,x},
\quad 2\le q<\infty,
\]
and for $q=\infty$ the same final lower bound follows by taking
any time in $[x^2,2x^2]$. Integration over $[2,\infty)$ yields, for every
$2\le q\le\infty$,
\[
 \|G_q^{H^{(1)}}(\eta)\|_{L_p(\R)}
 \ge\frac{e^{-1/4}\sqrt{\log2}}{16\sqrt{2\pi}}
          2^{(1-p)/p}(p-1)^{-1/p}\ge\frac{e^{-1/4-1/e}\sqrt{\log2}}{64\sqrt\pi}p'>\frac{p'}{512}\|\eta\|_{L_p(\R)},
\]
which completes the proof.
\end{proof}

\begin{proof}[Proof of Theorem~\ref{thm:sharp} and
               Corollary~\ref{cor:failure}]
The upper estimate in Theorem~\ref{thm:sharp} follows from Theorem~\ref{thm:main}.
The constants in the lower estimate are uniform in $q$. For $2\le q<\infty$,
\[
\frac{\Gamma(q)^{1/q}}q\ge \frac{1}{q} \left(\int_{q/2}^{q}s^{q}e^{-s}\frac{\dd s}{s} \right)^{1/q}\ge\frac1{2e}.
\]
At $q=\infty$, the corresponding coefficient in \eqref{sharp:eq-shell} is $1/e$, so the same lower bound holds.
For $p\ge2$, Proposition~\ref{sharp:shell}, including its real scalar
version, gives an absolute multiple of $\sqrt p$, since
$p/(2\sqrt{p-1})\ge\sqrt p/2$. Also
$\max\{p',\sqrt p\}\le\sqrt2\sqrt p$ in this range. For $1<p\le2$,
we have $\max\{p',\sqrt p\}=p'$, and
Proposition~\ref{sharp:tail} gives the required bound already in
dimension one. Isometric scalar embeddings establish the same lower
order for every nonzero Hilbert space.
This proves the dimension-uniform order $\max\{p',\sqrt p\}$,
with lower comparison constant $1/512$, independent of $q$,
and the upper bounds are those in \eqref{eq:main}.

Now fix $2<q\le\infty$, with the convention $1/\infty=0$. For $p\ge2$, the lower bound
in Proposition~\ref{sharp:shell} implies
\[
 \frac{\sup_{d\ge1}\mathrm L_{c,q,p}^{H^{(d)}}(\R)}
      {\max\{p',p^{1/q}\}}
 \ge\frac1{16e}p^{1/2-1/q}\longrightarrow\infty,
 \quad p\to\infty.
\]
Here $\max\{p',p^{1/q}\}\le2p^{1/q}$ and the real scalar shell
bound is at least $\sqrt p/(8e)$.
Consequently, no constant depending only on $q$ can make the estimate
with order $\max\{p',p^{1/q}\}$ hold simultaneously for all $p$ and
all Euclidean dimensions. The same conclusion holds for complex
scalars and for a fixed real two-dimensional Hilbert space.
The proof is complete.
\end{proof}

The dimensions in the shell construction vary with $p$ and with the
accuracy of approximation. Accordingly, this argument makes no claim
about the large-$p$ order in any one fixed Euclidean dimension when
dimension-dependent constants are allowed.
For $q=2$, Theorem~\ref{thm:main} gives the dimension-free upper
bound for the heat semigroup. For the classical Poisson semigroup,
the dimension-free forward estimate has order
$\max\{p',p^{1/q}\}$ for $2\le q<\infty$; see
\cite[Corollary~1.6]{XuJEMS}.

\section{Applications}\label{sec:applications}

We illustrate Theorem~\ref{thm:main} by varying the underlying
measure space $\Omega$, the generator $A$, and the Hilbert
space $\mathcal{H}$.
In all the examples below, the numerical bounds in
\eqref{eq:main} apply for $1<p<\infty$ and $2\le q\le\infty$.
The projection $F$ is the fixed point projection of the
corresponding semigroup.
The standard form constructions can be found in
\cite{Davies,FOT}.

\subsection{Underlying measure spaces}\label{sec7.1}

The following spaces carry natural symmetric submarkovian
heat semigroups.

\begin{itemize}

\item
Let $\Omega=M$ be a smooth complete Riemannian manifold
without boundary, with Riemannian measure $\mu$, and take
$A=-\Delta_M$ and $T_t=e^{t\Delta_M}$, the minimal heat
semigroup; see \cite{Davies,Strichartz}.
This includes $\R^d$, the unit sphere $\mathbb S^d$,
the flat torus $\mathbb T^d$, hyperbolic space $\mathbb H^d$,
complete Einstein manifolds, Calabi--Yau manifolds with a
complete Ricci-flat K\"ahler metric, Cartan--Hadamard
manifolds, and complete manifolds with finitely many ends.
No curvature bound or stochastic completeness assumption
is needed.
More generally, let $M$ be as above and let $V$ be a
smooth real-valued function on $M$.
Equip $M$ with the weighted measure $ \dd\mu=e^{-V}\dd\mathrm{vol}$,
and take the minimal selfadjoint realization of the
weighted Laplacian
\[
 A=-e^V\operatorname{div}(e^{-V}\nabla)
   =-\Delta_M+\langle\nabla V,\nabla\rangle.
\]
Its heat semigroup is symmetric and submarkovian and
has a measurable heat kernel; see \cite{Davies,FOT}.
The choice $V=0$ recovers the preceding example.
For $M=\R^d$ and $V(x)=|x|^2/2$, normalizing $\mu$
gives standard Gaussian measure, and
$A=-\Delta+x\cdot\nabla$ is the Ornstein--Uhlenbeck
operator; see \cite[Example~A.9]{XuJEMS}.

\item
Let $(\Omega,d,\mu)$ be an $\mathrm{RCD}^{*}(K,N)$ space,
where $K\in\R$ and $1\le N<\infty$.
Take $A=-\Delta$, the nonnegative generator of the
canonical heat flow associated with the Cheeger energy.
Its heat semigroup has a symmetric positive Markov
kernel; see \cite{JLZ}.

\item
Let $\Omega=E$ be an infinite-dimensional separable real
Banach space equipped with a nondegenerate centered
Gaussian probability measure $\gamma$.
The classical Ornstein--Uhlenbeck semigroup is
\[
 T_t(f)(x)=
 \int_E f(e^{-t}x+\sqrt{1-e^{-2t}}\,y)
 \dd\gamma(y).
\]
It is symmetric and Markovian and has a measurable
transition kernel; see \cite[Section~3]{LunardiPallara}.
Taking $A$ to be its nonnegative selfadjoint generator
gives an infinite-dimensional Gaussian example.
Here $E$ is the underlying space, while the functions
still take values in $\mathcal{H}$.

\item
Let $\Omega=G$ be a connected, simply connected stratified
nilpotent Lie group with Haar measure.
For horizontal left-invariant vector fields
$X_1,\ldots,X_m$ forming a basis of the first stratum, take
the standard selfadjoint sub-Laplacian $ A=-\sum_{j=1}^m X_j^2$.
This includes Heisenberg groups and groups of Heisenberg
type; see \cite{VSC}.

\item
The classical weighted polynomial heat semigroups on
the unit ball $\mathbb{B}^d= \{x\in\R^d:|x|<1\}$ and the standard simplex $\Sigma_d= \{x\in\R^d:  x_1+\dots+x_d<1,\ x_1,\dots,x_d>0\}$ are included.
Here the reference measures are the classical Jacobi
weights, and $A$ is the canonical selfadjoint polynomial
generator described in \cite{KPP}.

\item
Let $\Omega=V$ be a countable weighted space and let $P$
be a reversible Markov transition operator.
Then $A=\mathrm{Id}_{V}-P$ generates the symmetric Markov semigroup $ T_t=e^{-t(\mathrm{Id}_{V}-P)}$. This includes continuous-time reversible random walks
on weighted graphs and, after a rescaling of time,
the standard discrete heat semigroup on $\mathbb Z^d$; see \cite{CG1998} for the random walks.

\item
Let $\Omega$ be a finite compact metric graph equipped
with arc-length measure.
Take $A=-\dd^2/\dd s^2$ on each edge, with continuity
and Kirchhoff conditions at the vertices.
The resulting selfadjoint Laplacian generates a
symmetric Markov heat semigroup;
see \cite[Section~6.5]{Mugnolo}.

\item
The standard Brownian semigroup on the Sierpi{\'n}ski
gasket, equipped with its natural Hausdorff measure,
provides a fractal example.
Its generator is the corresponding fractal Laplacian;
see~\cite{BarlowPerkins}.

\end{itemize}

\subsection{Symmetric submarkovian generators}
\label{sec7.2}

The theorem also applies to the following classical
generators, with the indicated reference measures
and selfadjoint realizations.

\begin{itemize}

\item
Let $w$ be a Muckenhoupt $A_2$ weight on $\R^d$, and
let $a(x)$ be a real symmetric measurable matrix satisfying
\[
 \lambda w(x)|\xi|^2
 \le \langle a(x)\xi,\xi\rangle
 \le \Lambda w(x)|\xi|^2, \quad \xi\in\R^d,
\]
almost everywhere, where $0<\lambda\le\Lambda<\infty$.
On $(\R^d,w\dd x)$, take $A=-w^{-1}\operatorname{div}(a\nabla)$,
defined by the closure of the form
\[\int_{\R^d}\langle a\nabla f,\nabla f\rangle\dd x,\quad f\in C_c^\infty(\R^d).\]
These are degenerate elliptic operators of
Fabes--Kenig--Serapioni type; see \cite{FKS}.
Their heat semigroups are symmetric and submarkovian
and admit measurable heat kernels; see \cite{AtaeiNystrom}.
The case $w=1$ includes real symmetric uniformly elliptic
operators with bounded measurable coefficients.

\item
Let $A_0$ satisfy the assumptions of
Theorem~\ref{thm:main}, and let $\phi$ be a Bernstein
function.
Then $A=\phi(A_0)$, defined by spectral calculus,
generates a symmetric submarkovian kernel semigroup
by Bochner subordination; see \cite{SSV}.
In particular, this includes $A=A_0^\alpha$,
$0<\alpha<1$.
For $A_0=-\Delta$ on $\R^d$, these are fractional
Laplacians, with $\alpha=1/2$ giving the classical
Poisson semigroup.

\item
Let $\Omega\subset\R^d$ be open and let $V\ge0$ be locally
integrable.
The Schr\"odinger operator $A=-\Delta_\Omega+V$ is defined
by the closed form
\[
 \int_\Omega(|\nabla f(x)|^2+V(x)|f(x)|^2)\dd x,
 \quad
 f\in W_0^{1,2}(\Omega)\cap L_2(\Omega,V\dd x).
\]
Its heat semigroup is symmetric and submarkovian;
see \cite{Davies}.

\item
Let $R$ be a finite reduced root system in $\R^d$,
$R_+$ a positive subsystem, and $\kappa\ge0$ an invariant
multiplicity.
Writing $\sigma_\alpha$ for reflection across
$\alpha^\perp$, define
\[
 D_j(f)(x)=\partial_jf(x)
 +\sum_{\alpha\in R_+}\kappa(\alpha)\alpha_j
   \frac{f(x)-f(\sigma_\alpha x)}{\langle\alpha,x\rangle}.
\]
The Dunkl Laplacian
\[A=-\Delta_\kappa=-\sum_{j=1}^dD_j^2\]
has its standard nonnegative selfadjoint realization on
$L_2(\R^d,\dd\mu_\kappa)$, where
\[
 \dd\mu_\kappa(x)=
 \prod_{\alpha\in R_+}
 |\langle\alpha,x\rangle|^{2\kappa(\alpha)}\dd x.
\]
Its heat semigroup is symmetric and Markovian;
see \cite[Sections~2.2 and~4.2]{Rosler}.

\item
The Hermite operator $A=-\Delta+|x|^2$ on $\R^d$ with Lebesgue measure generates a symmetric submarkovian semigroup;
its heat kernel is given by Mehler's formula; see \cite{Thangavelu,XuJEMS}.

\item
The standard Grushin operator $ A=-\Delta_x-|x|^2\partial_u^2$
on $\R^d\times\R$ with Lebesgue measure generates
a symmetric Markov heat semigroup;
see \cite{DziubanskiJotsaroop}.

\item
For $\alpha>-1$, the Bessel operator
\[
 A=-\frac{\dd^2}{\dd x^2}
   -\frac{2\alpha+1}{x}\frac{\dd}{\dd x}
\]
on $(0,\infty)$ with measure $x^{2\alpha+1}\dd x$
is taken in the selfadjoint realization defined by the modified Hankel transform.
Its heat semigroup is symmetric and Markovian and admits a
measurable heat kernel. See \cite[Proposition~6.2 and (6.2)]{NowakStempak}.

\item
For $\alpha>-1$, the Laguerre polynomial operator
\[
 A=-x\frac{\dd^2}{\dd x^2}
   -(\alpha+1-x)\frac{\dd}{\dd x}
\]
on $(0,\infty)$ with measure $x^\alpha e^{-x}\dd x$
is taken in its polynomial selfadjoint realization;
see \cite{Dinger}.
The Laguerre operator of Hermite type
\[
 A=\frac12\left(
   -\frac{\dd^2}{\dd x^2}
   +x^2+\frac{\alpha^2-1/4}{x^2}
 \right)
\]
with Lebesgue measure is also included for
$\alpha\ge1/2$ with the Friedrichs realization, and
for $\alpha=-1/2$ with the Neumann realization;
see \cite{NowakStempak,XuJEMS}.

\item
For $\alpha,\beta>-1$, take $\Omega=[-1,1]$ with $ \dd\mu(x)=(1-x)^\alpha(1+x)^\beta\dd x$.
Writing $w(x)=(1-x)^\alpha(1+x)^\beta$, let $A$ be
the selfadjoint closure on $L_2(\Omega,\mu)$ of
\[
 A(f)=-\frac1w\frac{\dd}{\dd x}
          \left((1-x^2)w\frac{\dd f}{\dd x}\right),
\]
initially defined on algebraic polynomials.
Its heat semigroup is symmetric and Markovian;
see \cite[Section~5]{KPP}.\footnote{\ Tensor products of the preceding Bessel, Laguerre,
and Jacobi semigroups give the corresponding multidimensional examples.}

\item
On the half-space
$\R^d_+=\{x\in\R^d:x_d>0\}$ with Lebesgue measure,
both the Neumann Laplacian and the Dirichlet Laplacian
generate symmetric submarkovian heat semigroups.
Their form domains are $W^{1,2}(\R^d_+)$ and
$W_0^{1,2}(\R^d_+)$, respectively; see \cite{Davies}.
\end{itemize}

\subsection{Hilbert-valued functions}\label{sec7.3}

Every preceding example allows any Hilbert space $\mathcal{H}$.
The semigroup acts on the variable in $\Omega$ through
$T_t\otimes\mathrm{Id}_{\mathcal{H}}$, with the same numerical bounds.

\begin{itemize}

\item
Taking $\mathcal{H}=\R^m$ or $\C^m$ with the Euclidean norm gives
vector-valued estimates, uniformly in $m$.
Taking $\mathcal{H}=\ell_2$ gives sequence-valued estimates.
Hilbert direct sums are included with the same constants.

\item
One may take $\mathcal{H}=L_2(Y,\nu)$ for an auxiliary
$\sigma$-finite measure space $(Y,\nu)$.
The estimates then hold on $L_p(\Omega;L_2(Y,\nu))$.
Other choices include $\mathcal{H}=W^{1,2}(D)$ and
$\mathcal{H}=H^s(\R^m)$, with their usual Hilbert norms.
The operator $A$ acts only on the $\Omega$ variable.

\item
Let $\mathcal N$ be a semifinite von Neumann algebra
with a normal semifinite faithful trace $\tau$.
The choice
\[
 \mathcal{H}=L_2(\mathcal N,\tau),
 \quad \|a\|_{\mathcal{H}}=\tau(a^*a)^{1/2},
\]
gives estimates on the Bochner space
$L_p(\Omega;L_2(\mathcal N,\tau))$;
see \cite{PisierXuLp}.
In particular, $\mathcal N=M_m(\C)$ gives
Hilbert--Schmidt matrix-valued functions, uniformly
in $m$, and the infinite-dimensional class $S_2$
is included as well.
The underlying semigroup still acts on the classical
space $\Omega$.

\end{itemize}

\subsection{A further problem}\label{sec7.4}

The applications in
Subsections~\ref{sec7.1}--\ref{sec7.3}
suggest the following problem in the setting of
\cite{XuJEMS}.

\begin{problem}\label{prob:cotype}
Let $X$ be a Banach space of martingale cotype $q_0 \ge 2$.
Can one find a constant $C_{X,q_0}$ such that, for every semigroup satisfying the assumptions of Theorem~\ref{thm:main},
\[
\|G_q^T(f)\|_{L_p(\Omega)} \le C_{X,q_0} \max\{p', \sqrt{p}\} \|f - F(f)\|_{L_p(\Omega; X)}
\]
holds for all $1<p<\infty$, all $q_0\le q\le\infty$, and all $f\in L_p(\Omega;X)$?
Here $G_q^T$ is defined using the norm of $X$, and $C_{X,q_0}$ is independent of $p,q$, the underlying space $\Omega$, and the semigroup $(T_t)_{t>0}$.
\end{problem}

\smallskip

{\footnotesize{\noindent \textbf{Acknowledgments.}
The authors acknowledge the use of GPT-6.1 in the preparation of this manuscript and assume full responsibility for its content.
Xudong Lai is supported by NNSF of China (12322107, 12271124, W2441002) 
and Heilongjiang NSF of China (YQ2022A005).
Bo Li is supported by NNSF of China (12471094, 12571121), Zhejiang NSF of China (LMS26A010015).
Ji Li is supported by ARC (DP260100485).
}}

\begingroup
\footnotesize
\setlength{\parskip}{0pt}
\setlength{\itemsep}{0pt}
\setlength{\parsep}{0pt}
\setlength{\topsep}{0pt}
\setlength{\partopsep}{0pt}
\setlength{\labelsep}{0.3em}
\renewcommand{\baselinestretch}{0.95}\selectfont
\emergencystretch=1em

\endgroup

\footnotesize{
\noindent  Xudong Lai (xudonglai@hit.edu.cn, xudonglai@mail.bnu.edu.cn), Institute for Advanced Study in Mathematics, Harbin Institute of Technology, Harbin 150001, China;
Zhengzhou Advanced Research Institute, Harbin Institute of Technology, Zhengzhou 450000, China.}

\smallskip

\noindent Bo Li (bli@zjxu.edu.cn), Institute of Mathematics, Jiaxing University, Jiaxing 314001, China.

\smallskip

\noindent Ji Li (ji.li@mq.edu.au), School of Mathematical and Physical Sciences, Macquarie University, NSW 2109, Australia.

\end{document}